\RequirePackage{fix-cm}
\documentclass{svjour3}                     % onecolumn (standard format)
\usepackage{graphicx}
\usepackage[a4paper,hscale=0.7,vscale=0.75,centering]{geometry}
\usepackage{amssymb,mathrsfs,esint}
\usepackage{color}
\usepackage{cite}

\renewcommand{\r}{\mathbb{R}}
\newcommand{\rd}{\mathbb{R}^d}

\newcommand{\E}{\mathbb{E}}

\renewcommand{\P}{\mathbb{P}}
\renewcommand{\tilde}{\widetilde}
\renewcommand{\epsilon}{\varepsilon}

\newtheorem{prop}{Proposition}
\newtheorem{thm}[prop]{Theorem}
\newtheorem{coroll}[prop]{Corollary}
\journalname{PTRF}
\begin{document}

\title{Reflection couplings and contraction rates for diffusions}


\titlerunning{Reflection couplings and contraction rates}        % if too long for running head

\author{Andreas Eberle %etc.
}


\institute{Andreas Eberle \at
              Institute for Applied Mathematics, University of Bonn, Endenicher Allee 60, 53115 Bonn, Germany \\
              \email{eberle@uni-bonn.de}           %  \\
}

\date{Received: date / Accepted: date}

\maketitle

\begin{abstract}
We consider contractivity for diffusion semigroups w.r.t.\ Kantoro\-vich ($L^1$ Wasserstein) distances
based on appropriately chosen concave functions. These distances are inbetween total variation and
usual Wasserstein distances. It is shown that by appropriate explicit choices of the underlying
distance, contractivity with rates of close to optimal order can be obtained in several fundamental
classes of examples where contractivity w.r.t.\ standard Wasserstein distances fails. Applications
include overdamped Langevin diffusions with locally non-convex potentials, products of these
processes, and systems of weakly interacting diffusions,
both of mean-field and nearest neighbour type.
\end{abstract}

\section{Introduction}

Consider a diffusion process $(X_t)_{t\geq 0}$ in $\mathbb R^d$ defined by a stochastic differential equation
\begin{equation}
\label{eq:1}
dX_t\:=\:b(X_t)\:dt+\sigma\:dB_t.
\end{equation}
Here $(B_t)_{t\geq 0}$ is a $d$-dimensional Brownian motion, $\sigma\in \mathbb R^{d\times d}$ is a constant $d\times d$ matrix with $\det \sigma>0$, and $b:\mathbb R^d\to\mathbb R^d$ is a locally Lipschitz continuous function.  We assume that the unique strong solution of (\ref{eq:1}) is non-explosive
for any initial condition, which is essentially a consequence of the assumptions imposed further below. The transition kernels of the diffusion process on $\mathbb R^d$ defined by (\ref{eq:1}) will be denoted by $p_t(x,dy)$.\medskip

Contraction properties of the transition semigroup $(p_t)_{t\ge 0}$
%of the diffusion process $(X_t)_{t\ge 0}$
have been studied by various approaches. In particular,
$L^2$ and entropy methods (e.g.\ spectral gap estimates, logarithmic Sobolev and transportation inequalities) yield bounds that both are relatively
stable under perturbations and applicable in high dimensions, cf.\ e.g.\ \cite{ABC, BE, B, BCG, BGL, BGG, Roy, Wang}. On the other hand,
coupling methods provide a more intuitive probabilistic
understanding of convergence to equilibrium
\cite{LR, Lindvall, Thorisson, CL, Chen, CW, Wang, HM, HMS}.
In contrast to $L^2$ and entropy methods, bounds resulting
from coupling methods typically hold for arbitrary initial
values $x_0\in\mathbb{R}^d$.
In many applications, couplings are used to bound the total
variation distances $d_{TV}(\mu p_t,\nu p_t)$ between the
laws $\mu p_t$ and $\nu p_t$ of $X_t$ w.r.t.\ two different initial distributions
$\mu$ and $\nu$ at a given time
$t\ge 0$ , cf.\ \cite{Lindvall, LR}. Typically,
however, the total variation distance is decaying substantially only after a certain amount of time. This is also manifested in cut-off
phenomena \cite{Diaconis, Peres, DSC, ChenSC}.\medskip

Alternatively, it is well-known that synchronuous couplings
(i.e., couplings given by the flow of the s.d.e.\
(\ref{eq:1})) can be used to show that the map
$\mu\mapsto\mu p_t$ is exponentially contractive w.r.t.\
$L^p$ Wasserstein distances $W^p$ for any $p\in [1,\infty )$ if, for example,
$(X_t)$ is an overdamped Langevin diffusion with a strictly
convex potential $U\in C^2(\mathbb{R}^d)$, i.e.,
$\sigma =I_d$ and $b=-\nabla U/2$, see e.g.\ \cite{BGG}. This leads to an elegant and powerful approach to convergence to equilibrium
and to many related results if applicable. However, it
has been pointed out in \cite{RenesseSturm} that strict
convexity of $U$ is also a necessary condition for
exponential contractivity w.r.t.\ $W^p$. This seems
to limit the applicability substantially.\medskip

Here, we are instead considering exponential contractivity
w.r.t.\ Kantorovich ($L^1$ Wasserstein) distances $W_f$
based on underlying distance functions of the form
$$d_f(x,y)\ =\ f(\| x-y\| )\qquad\mbox{ on }\ \mathbb{R}^d,$$
and, more generally,
$$d_f(x,y)\ =\ \sum_{i=1}^nf_i(\| x^i-y^i\| )
\qquad\mbox{ on }\ \mathbb{R}^{d_1}\times\cdots\times\mathbb{R}^{d_n} ,$$
where $f,f_i:[0,\infty )\to [0,\infty )$ are strictly
increasing {\em concave} functions, cf.\ Sections \ref{subsectionRC}
and \ref{subsectioncrc} below for details. For proving
exponential contractivity, we will apply a reflection
coupling on $\rd$ and an (approximate) componentwise
reflection coupling on products of Euclidean spaces.
It will become clear by the proofs below, that for
distances based on concave functions $f,f_i$, these
couplings are superior to synchronuous couplings, whereas the
synchronuous couplings are superior w.r.t.\ the Wasserstein
distances $W^p$ for $p>1$, cf.\ e.g.\ Lemma \ref{lem3A}.\medskip

The idea to study contraction properties w.r.t.\ Kantorovich distances based on
concave distance functions appears in Chen and Wang \cite{WangNeumann, CWa, CW} and Hairer and Mattingly \cite{HM}. In \cite{CW},
similar methods are applied to estimate
spectral gaps of diffusion generators on $\mathbb{R}^d$ and on
manifolds. In \cite{HM} and \cite{HMS},
Hairer, Mattingly and Scheutzow apply Wasserstein
distances based on particular concave distance functions
to prove exponential ergodicity in infinite dimensional situations.
The key idea below is to obtain more quantitative results by ``almost'' optimizing the choice of the functions $f$ and $f_i$ to obtain large contraction rates. In the case $n=1$, this idea has also
been exploited in
\cite{CW} to derive lower bounds for spectral gaps.
The novelty here is that we suggest a simple and very {\em explicit}
choice for $f$ that leads to close to optimal results
in several examples. Furthermore, by a new extension to the product case
based on an approximate componentwise reflection coupling,
we obtain dimension free contraction
results in product models and perturbations thereof
without relying on convexity.\medskip

Before stating the general results, we consider some
examples illustrating the scope of the approach:

%We are interested in upper bounds for Kantorovich-Rubinstein-Wasserstein distances of the distributions $\mu p_t$ and $\nu p_t$ at a given time $t\geq 0$ w.r.t.\ two different initial distributions $\mu$ and $\nu$.\smallskip

\begin{example}[{\bf Overdamped Langevin dynamics with locally non-convex potential}]\label{exlangevin}
Suppose that $\sigma=I_d$ and $b(x)=-\frac 12\nabla U(x)$ for a function
$U\in C^2(\mathbb R^d)$ that is
%(i.e. $\nabla^2 U\geq K \cdot I_d$ for some $K>0$)
{\em strictly convex outside a given ball} $B\subset
\mathbb R^d$. Then $Z:=\int \exp(-U(x))dx$ is finite, and
the probability measure $$d\mu
\ =\ Z^{-1}\exp(-U)\,dx$$ is a stationary distribution for the diffusion
process $(X_t)$. Corollary \ref{cor16} below yields exponential contractivity for the transition semigroup $(p_t)$ with an explicit rate w.r.t.\ an appropriate
Kantorovich distance $W_f$. As a consequence, we
obtain dimension-independent upper bounds for the standard $L^1$
Wasserstein distances between the laws $\nu p_t$ of $X_t$ and $\mu$
for arbitrary initial distributions $\nu$ and $t\geq 0$.
These bounds are of of optimal order in $R,L\in [0,\infty )$
and $K\in (0,\infty )$ if $(x-y)\cdot(\nabla U(x)-\nabla U(y))$ is bounded from below by
$-L|x-y|^2$ for $|x-y|<R$ and by $K|x-y|^2$ for $|x-y|\geq R$.
\end{example}

\begin{example}[{\bf Product models}]\label{exproduct}
For a diffusion process $X_t=(X_t^1,\ldots ,X_t^n)$ in
$\mathbb R^{n\cdot d}$ with independent Langevin diffusions
$X^1,\ldots ,X^n$ as in Example \ref{exlangevin},
Theorem \ref{thmA} below yields exponential contractivity in an
appropriate Kantorovich distance with rate
$c=\min (c_1,\ldots ,c_n)$ where $c_1,\ldots ,c_n$ are the
lower bounds obtained for the contraction rates of the
components.
\end{example}

\begin{example}[{\bf Systems of interacting diffusions}]
\label{exinteracting}
More generally, consider a system
$$dX_t^i\ =\ -\frac 12 \nabla U(X_t^i)\, dt\, -\,
\frac{\alpha}{n}\sum_{j=1}^n\nabla V(X_t^i-X_t^j)\, dt
\, +\, dB_t^i, \ \ \ i=1,\ldots ,n,$$
of $n$ interacting
diffusion processes in $\mathbb R^d$ where $U\in C^2(\mathbb R^d)$
is strictly convex outside a ball, $V\in C^2(\mathbb R^d)$
has bounded second derivatives, and $B^1,\ldots ,B^n$ are
independent Brownian motions in $\mathbb R^d$. Then
Corollary \ref{corC} below shows that for $\alpha$
sufficiently small, exponential contractivity holds in
an appropriate Kantorovich distance with a rate that does
not depend on $n$.
\end{example}

We now introduce briefly the couplings to be considered in the proofs below:

%\subsection{Couplings for multidimensional diffusions}
A \emph{coupling by reflection} of two solutions of (\ref{eq:1}) with initial distributions $\mu$ and $\nu$ is a diffusion process $(X_t,Y_t)$ with values in $\mathbb R^{2d}$ defined by $(X_0,Y_0)\sim \eta$ where $\eta$ is a coupling of $\mu$ and $\nu$,
\begin{eqnarray}
dX_t&=&b(X_t)\:dt+\sigma\:dB_t\qquad\qquad\qquad\mbox{for }t\geq 0,\\
dY_t&=&b(Y_t)\:dt+\sigma (I-2e_te_t^\top)\:dB_t\quad\mbox{ for $t<T$,}\quad Y_t\ =\ X_t\mbox{ for }t\geq T.\label{eq:2}
\end{eqnarray}
Here $e_te_t^\top$ is the orthogonal projection onto the unit vector
$$e_t:=\sigma^{-1}(X_t-Y_t)/|\sigma^{-1}(X_t-Y_t)|,$$
and
$T=\inf\{t\geq 0\, :\, X_t=Y_t\} $ is the coupling time, i.e., the
first hitting time of the diagonal $\Delta=\{(x,y)\in\mathbb
R^{2d}:x=y\}$, cf.~\cite{LR,CL}. The reflection coupling can be
realized as a diffusion process in $\mathbb R^{2d}$, and the
marginal processes $(X_t)_{t\geq 0}$ and $(Y_t)_{t\geq 0}$ are
solutions of (\ref{eq:1}) w.r.t.\ the Brownian motions $B_t$  and
$$\check B_t=\int_0^t(I_d-2\mathbb{I}_{\{s<T\}}e_se_s^\top)\:dB_s.$$
Notice that by L\'evy's characterization, $\check B$ is
indeed a Brownian motion since the process $I_d-2\mathbb{I}_{\{s<T\}}e_se_s^\top$
takes values in the orthogonal matrices.
The
difference vector $$Z_t\ :=\ X_t-Y_t$$ solves the s.d.e.
\begin{eqnarray}
\label{eq:3}
dZ_t &=&(b(X_t)-b(Y_t))\:dt+{2}{|\sigma^{-1}Z_t|^{-1}}Z_t\:dW_t\quad\mbox{ for }t<T,\\ Z_t&=&0\quad\mbox{ for }t\ge
T,\nonumber
\end{eqnarray}
w.r.t.\ the \emph{one-dimensional} Brownian motion
$$W_t\:=\:\int_0^te_s^\top\:dB_s.$$
%\medskip\\

A {\em synchronuous coupling} of two solutions of (\ref{eq:1})
is defined correspondingly with $e_t\equiv 0$, i.e.,
the same noise is applied both to $X_t$ and $Y_t$. Below we will also consider {\em mixed couplings} that are reflection
couplings for certain values of $Z_t$, synchronuous couplings
for other values of $Z_t$, and mixtures of both types of couplings for $Z_t$ in an intermediate region.
Notice that the standard reflection coupling introduced
above is a synchronuous coupling for $t\ge T$, i.e.,
if $Z_t=0$~! \medskip

More generally, we will consider
couplings for diffusion processes on product spaces
(such as in Examples \ref{exproduct} and \ref{exinteracting}) that are approximately
{\em componentwise reflection couplings}, i.e., the $i$-th
component $(X^i_t,Y^i_t)$ of the coupling $(X_t,Y_t)$ is defined
similarly to (\ref{eq:2}) provided $|X_t^i-Y_t^i|\ge\delta $ for a given
 constant $\delta >0$, cf.\ Section \ref{secprod} below.\medskip

For diffusion processes with non-constant diffusion matrix $\sigma (x)$, the reflection coupling
should be replaced by the {\em Kendall-Cranston coupling} w.r.t.\ the intrinsic Riemannian metric
$G(x)=\left(\sigma (x)\sigma (x)^T\right)^{-1}$ induced by the diffusion coefficients, cf.\
\cite{Kendall, Cranston, Hsu, Wang}. Here, we restrict ourselves to the case of constant diffusion
matrices where the Kendall-Cranston coupling coincides with the standard coupling by
reflection.\medskip

The main results of this paper are stated in Section 2 for
reflection coupling, and in Section 3 for componentwise
reflection coupling on product spaces. The
proofs are contained in Sections 4, 5 and 6. A part of the results in Section \ref{secmain1} have been announced in the Comptes Rendus Note \cite{EberleCR}.

\section{Main results for reflection coupling}\label{secmain1}

\subsection{Reflection couplings and contractivity on $\mathbb R^d$}\label{subsectionRC}

Lindvall and Rogers \cite{LR} introduced coupling by reflection in order to derive upper bounds for the total variation distance of the distributions of $X_t$ and $Y_t$ at a given time $t\geq 0$. Here we are instead considering the Kantorovich-Rubinstein ($L^1$-Wasserstein) distances
\begin{equation}
\label{eq:4} W_f(\mu,\nu)\:=\:\inf\limits_\eta\int
d_f(x,y)\,\eta(dx\:dy),\,\quad d_f(x,y)\:=\:f(\|x-y\|)\quad
(x,y\in\mathbb R^d),\end{equation} of probability measures $\mu,
\nu$ on $\mathbb R^d$, where the infimum is over all couplings
$\eta$ of $\mu$ and $\nu$, $f:[0,\infty)\to[0,\infty)$ is an
appropriately chosen concave increasing function with $f(0)=0$, and
$\|z\|=\sqrt{z\cdot Gz}$ with $G\in\mathbb R^{d\times d}$ symmetric
and strictly positive definite. Typical choices for the norm are the
Euclidean norm $\|z\|=|z|$ and the intrinsic metric
$\|z\|=|\sigma^{-1}z|$ corresponding to $G=I_d$ and
$G=(\sigma\sigma^\top)^{-1}$ respectively.

\begin{remark}[{\bf Interpolating between total variation and Wasserstein distances}]
For the choice of the function $f$ there are two extreme cases with minimal and maximal concavity:
\begin{enumerate}
\item Choosing $f(x)=x$ yields the {\em standard Kantorovich
($L^1$~Wasserstein) distance} $W_f=W^1$. In this case it is well known that if, for example, $G=\sigma =I_d$ and $b(x)=-\nabla U(x)/2$, then the transition kernels $p_t(x,dy)$ of the diffusion process $(X_t)$ satisfy $$W_f(\mu p_t,\nu p_t)\ \leq\  e^{-K t/2}\, W_f(\mu,\nu)\qquad
\mbox{for any $\mu,\nu$ and }t\geq 0,$$ provided $\nabla^2U\geq K\cdot I_d$ holds globally. This condition is also sharp in the sense that if $U$ is not globally strictly convex, then contractivity of $p_t$ w.r.t.\ $W_f$ does not hold, cf.\ Sturm and von
Renesse \cite{RenesseSturm}.
\item On the other hand,
choosing $f(x)=\mathbb{I}_{(0,\infty )}(x)$ yields the {\em total variation distance} $W_f=d_{TV}$. In this case,
   $$W_f(\mu p_t,\nu p_t)\ \leq\ \P[T>t]\qquad
   \mbox{for any $\mu,\nu$ and }t\geq 0,$$
   but there is no strict contractivity of $p_t$ w.r.t.\ $d_{TV}$ in general. Indeed, in many applications $d_{TV}(\mu p_t,\nu p_t)$ only decreases substantially after a certain amount
of time (``cut-off phenomenon'').
\end{enumerate}
\end{remark}

By choosing for $f$ an appropriate concave function, exponential contractivity w.r.t.\ $W_f$ may hold even without global convexity, cf.\ \cite{CW}. We now explain how the function $f$
can be chosen in a very explicit way such that the
obtained exponential decay rate w.r.t.\ the Kantorovich
distance $W_f$ differs from the maximal decay rate that
we can achieve by our approach based on reflection coupling
only by a constant factor.\medskip

%The idea to study Wasserstein contractivity w.r.t.\
%concave distance functions goes back to Chen and Wang \cite{CW},
%where it is implicitly contained in the proofs. Indeed, in \cite{CW}
%and \cite{W}, Chen and Wang apply very similar methods to estimate
%spectral gaps of diffusion generators on $\mathbb{R}^d$ and on
%manifolds. Related arguments have also been applied in \cite{HM} to
%quantify exponential ergodicity in infinite dimensional situations.
%
%The function $f$ will then be chosen in an
%The key idea in the result below is to
%``almost'' optimal way to obtain a large exponential decay rate. This idea has
%also been exploited intensively by Chen and Wang in
%\cite{CW} to derive lower bounds for spectral gaps.
%The novelty here is that we suggest a simple and very {\em explicit}
%choice for $f$ that leads to close to optimal results
%in several examples.
At first, similarly to Lindvall and Rogers \cite{LR}, let us define for $r\in(0,\infty)$:
\[\kappa(r)\:=\:\inf\left\{ -2\,\frac{|\sigma^{-1}(x-y)|^2}{\|x-y\|^2}\,\frac{(x-y)\cdot G(b(x)-b(y))}{\|x-y\|^2}\, :\,
x,y\in\mathbb R^d\mbox{ s.t.\ }\|x-y\|=r\right\},\]
i.e., $\kappa (r)$ is the largest constant such that
\begin{equation}
\label{eqDEFK}
(x-y)\cdot G(b(x)-b(y))\ \le\ -\frac 12\kappa (r) \|x-y\|^4
/|\sigma^{-1}(x-y)|^2
\end{equation}
holds for any $x,y\in\mathbb R^d$ with $\|x-y\|=r$.
 Notice that if $\|\,\cdot\,\|$
is the intrinsic metric then the
factor $|\sigma^{-1}(x-y)|^2/\|x-y\|^2$ equals 1 . In Example \ref{exlangevin} with $G=I_d$, we have $$\kappa
(r)=\inf\left\{\int\nolimits_0^1\partial^2_{(x-y)/|x-y|}U((1-t)x+ty)\,
dt :x,y\in\mathbb R^d \mbox{ s.t.}\, |x-y| =r\right\}.$$
We assume
from now on that $\kappa (r)$ is a continuous function on $(0,\infty )$ satisfying
\begin{equation}\label{eqassk}
\liminf_{r\to\infty}\kappa(r)>0\mbox{\qquad and
\qquad}\int_0^1r\kappa (r)^-\, dr<\infty .
\end{equation}
In Example \ref{exlangevin} with $G=I_d$, this assumption
is satisfied if $U$ is strictly convex outside a ball.
\smallskip

Next, we define constants $R_0,R_1\in[0,\infty)$ with $R_0\leq R_1$ by
\begin{eqnarray}
\label{eqdefr0}
R_0& =& \inf\{R\geq 0\::\:\kappa(r)\geq 0\ \forall\,r\geq R\},\\
\label{eqdefr1}
R_1& =& \inf\{R\ge R_0\::\:\kappa(r)R(R-R_0)\geq 8\ \forall\, r\geq R\} ,
\end{eqnarray}
% \begin{remark}In the {overdamped Langevin} case with
%$G=\sigma=I_d$,
 %\begin{equation}\label{eq:5}
%\kappa(r)\:=\:\inf\left\{\int\limits_0^1\partial_{(x-y)/\|x-y\|}U((1-r)x+ry)dr\ :\ \mbox{ $x,y\in\mathbb R^d$ with $\|x-y\|%=r$}\right\}.
%\end{equation}
%In particular, the assumption $\,\liminf\limits_{r\to\infty}\kappa(r)>0$ is satisfied if $U$ is strictly convex outside a ball.
%\end{remark}
Notice that by (\ref{eqassk}), both constants are finite. We
now consider the particular distance function $d_f(x,y)=f(\|x-y\|)$ given by
\begin{eqnarray}
f(r)&=&\int\limits_0^r\varphi(s)g(s)\:ds,\qquad\qquad \mbox{where}\label{eq:8}
\\
\nonumber
\varphi(r)&=&\exp\left(-\frac 14\int\limits_0^rs\kappa(s)^-\:ds\right),\qquad
\Phi(r) \, =\, \int_0^r\varphi(s)\:ds,\\
g(r)&=& 1-\left.\frac 12\int\limits_0^{r\wedge R_1}\frac{\Phi(s)}{\varphi(s)}\:ds\right/\int\limits_0^{R_1}\frac{\Phi(s)}{\varphi(s)}\:ds.\nonumber
\end{eqnarray}
%Notice that $\Phi$ and $f$ are {\em concave}, because $\varphi$ and $g$ are
%{\em decreasing}. Moreover, $$\Phi(r)/2\leq f(r)\leq\Phi(r)\mbox{\qquad for any }r\geq 0.$$
%The choice of $f$ will become clear by the proof of the main result below. For the moment
Let us summarize some basic properties
of the functions $\varphi ,g$ and $f$:
\begin{itemize}
   \item $\varphi$ is decreasing, $\varphi(0)=1$, and $\varphi(r)=\varphi(R_0)$ for any $r\geq R_0$,
   \item $g$ is decreasing, $g(0)=1$, and $g(r)=\frac 12$ for any $r\geq R_1$,
   \item $f$ is concave, $f(0)=0$, $f'(0)=1$, and
   \begin{equation}
   \label{eq:9}
   \Phi(r)/2\:\leq\:f(r)\:\leq\:\Phi(r)\qquad\mbox{ for any }r\geq 0.
   \end{equation}
\end{itemize}
The last statement shows that
$d_f$ and $d_\Phi$ as well as $W_f$ and $W_\Phi$ differ at most by a factor $2$.\medskip

We will explain in Section \ref{secproofs1} below how the choice of $f$ is obtained
by trying to maximize the exponential decay rate. Let us now state our first main result which will be proven in
Section \ref{secproofs1}.
%Moreover, since
%\[\frac 12\Phi(R_0)\:\leq \:\frac 12\Phi(r)\:\leq\:f(r)\:\leq\:r,\]
%we have
%\[\frac 12\Phi(R_0)\cdot\|x-y\|\:\leq\:d_f(x,y)\:\leq\:\|x-y\|\qquad\mbox{ for any }x,y\in\mathbb R^d.\]
\medskip

\begin{thm}[{\bf Exponential contractivity of
reflection coupling}]\label{thmmain1}
Let $\alpha:=\sup\{|\sigma^{-1}z|^2\, :\,z\in\mathbb R^d\mbox{ with }\|z\|=1\}$, and define $c\in(0,\infty)$ by
\begin{equation}
\label{eq:10}
\frac 1c\ =\ \alpha\int\limits_0^{R_1}\Phi(s)\varphi(s)^{-1}\:ds\ =\
\alpha\int\limits_0^{R_1}\int\limits_0^s\exp\left(\frac 14\int\limits_t^su\kappa(u)^-\:du\right)\:dt\:ds\, .
\end{equation}
Then for the distance $d_f$ given by (\ref{eq:4}) and (\ref{eq:8}), the function
$t\mapsto e^{ct}\E[d_f(X_t,Y_t)]$
is decreasing on $[0,\infty)$.
\end{thm}\medskip
%\begin{remark}
%Note that $\alpha=1$ if $\|z\|=|\sigma^{-1}z|$ is the intrinsic metric of the diffusion.
%\end{remark}
 The theorem yields exponential contractivity at rate $c>0$ for the transition kernels $p_t$ of (\ref{eq:1}) w.r.t.\ the Kantorovich distance $W_f$. Moreover, it implies upper bounds for the standard Kantorovich ($L^1$~Wasserstein) distance $W^1=W_{\mbox{id}}$ w.r.t.\ the distance function $d(x,y)=\|x-y\|$:\smallskip

\begin{coroll}\label{cor16}
For any $t\geq 0$ and any probability measures $\mu,\nu$ on $\mathbb R^d$,
\begin{eqnarray}
W_f(\mu p_t,\nu p_t) &\leq & \exp ({-ct})\, W_f(\mu,\nu),\ \mbox{\qquad and }\label{eq:12}\\
W^1(\mu p_t,\nu p_t) & \leq & 2\varphi(R_0)^{-1}\exp ({-ct})\, W^1(\mu,\nu).\label{eq:12a}
\end{eqnarray}
\end{coroll}
Note that the second estimate follows from the first, since by the properties of $\varphi$ and $g$ stated above,
$\varphi (R_0)/2\le f^\prime \le 1$, and hence
%$\varphi (r)\ge \varphi (R_0)$ for any $r$ implies
\begin{equation}\label{eq113a}
\varphi(R_0)\|x-y\| /2\ \leq\  d_f(x,y)\ \leq\ \|x-y\| \qquad\mbox{ for any }x,y\in\mathbb R^d.
\end{equation}

The corollary yields an upper bound for mixing times
w.r.t.\ the Kantorovich distance $W^1$. For $\varepsilon >0$ let
\[\tau_{W^1}(\varepsilon)\ :=\ \inf\{t\geq 0\::\:W^1(\mu p_t,\nu p_t)\leq \varepsilon W^1(\mu,\nu)\ \forall\,\mu ,\nu\in
\mathcal M_1(\mathbb R^d )\} .$$
Then by Corollary \ref{cor16},
$$\tau_{W^1}(\varepsilon)\:\leq \: c^{-1}\log(2/(\varepsilon\varphi(R_0)))\quad\mbox{for any }\varepsilon>0.\]
The proofs of Theorem \ref{thmmain1} and Corollary \ref{cor16} are given in Section \ref{secproofs1} below.
% The bound can be used to control the $L^1-Wasserstein$ distance of the distribution at time $t$ and a stationary distribution %of the diffusion process.\end{remark}\medskip\\

\begin{remark}[{\bf Non-constant diffusion coefficients}]
The methods and results presented above have natural extensions to diffusion processes with smooth non-constant
diffusion matrices. In that case, one possibility is to
use an ad hoc coupling as in \cite{LR}, but this leads
to restrictive assumptions and bounds that are far from optimal. A better approach is
to switch to a Riemannian setup where the metric is the
intrinsic metric $G(x)=(\sigma (x)\sigma(x)^T)^{-1}$
given by the diffusion coefficients. The diffusion process $(X_t)$ can then be represented in the form
\begin{equation}
dX_t\ =\ \beta (X_t)\, dt\, +\, dB^G_t
\end{equation}
where $(B^G_t)$ is a Brownian motion on the Riemannian manifold $(\mathbb R^d,G)$, and $\beta$ is a modified
drift vector field. Now,
by replacing
the reflection coupling by the corresponding {Kendall-Cranston coupling} on $(\mathbb R^d,G)$, one can
expect similar results as above with $\kappa$ defined as
\[\kappa(r)\:=\:2r^{-1}\inf\left\{ -
\langle \gamma_{y,x}^\prime (r),{\beta (x)}\rangle
+\langle \gamma_{y,x}^\prime (0),{\beta (y)}\rangle
\,+\, \int_0^r{\rm Ric}(\gamma_{y,x}^\prime (s) ,\gamma_{y,x}^\prime (s) )\, ds   :
\|x-y\|=r\right\},\]
where $\gamma_{y,x}:[0,r]\to\mathbb R^d$ is the unit speed 
geodesic from $y$ to $x$ and ${\rm Ric}$ denotes the Ricci
curvature on $(\mathbb R^d,G)$, cf.\ \cite{Cranston, Wang}.
\end{remark}

\begin{remark}[{\bf Diffusions with reflection on smooth convex domains}]\label{rem:reflecteddiffusion}
The results above also apply to diffusion processes on a smooth bounded
domain $D\subseteq\mathbb R^d$ with normal reflection at the boundary
\cite{Tanaka, LionsSznitman, CranstonLeJan, BurdzyChenJones, Andres}.
In that case the SDE (\ref{eq:1}) is replaced by
\begin{equation}
dX_t\ =\ b(X_t)\, dt\, +\, n(X_t)\, d\ell_t\, +\,\sigma \, dB_t,
\end{equation}
where $n(x)$ is the interior normal vector at a boundary point $x$, and $(\ell_t)$ is the local time of $(X_t)$ on the boundary $\partial D$, i.e., $t\mapsto \ell_t$ is a non-decreasing process that increases only
at times when $X_t\in\partial D$. Consequently, in the Equation (\ref{eq:3})
for the coupling difference $Z_t=X_t-Y_t$, additional drift terms in the
directions $n(X_t)$ and $-n(Y_t)$ occur when one of the two copies is
at the boundary. Since for a convex domain, both $Z_t\cdot n(X_t)\le 0$
and $-Z_t\cdot n(Y_t)\le 0$, the reflection at the boundary improves the
upper bounds for $\| Z_t\|$ in the proofs below when choosing $G=I_d$.   
Therefore, the assertions of Theorem \ref{thmmain1} and Corollary \ref{cor16} hold true without further change
if we take the infimum in the definition of $\kappa$ 
only over $x,y\in D$ and choose $R_0$, $R_1$ respectively equal to the 
diameter of $D$ in case the infima in (\ref{eqdefr0}) or
(\ref{eqdefr1}) are over empty sets.
\end{remark}

\subsection{Consequences}\label{sub:cons}

We summarize some important consequences of exponential
contractivity w.r.t.\ Kantorovich distances as stated in
Corollary \ref{cor16}. These consequences are essentially well-known, cf.\ e.g.\ Joulin \cite{Joulin}, Joulin and Ollivier \cite{JO}, 
and Komorowski and Walczuk \cite{KW} for related results. For the
reader's convenience, the proofs are
nevertheless included in Section  \ref{secproofs1} below.
We assume that $\| z\| =|\sigma^{-1}z|$ is the
intrinsic metric, $b$ is in $C^1(\rd ,\rd )$, and
\begin{equation}
\label{eqFM}
\int |z|\, p_t(x_0,dz)\ <\ \infty
\end{equation}
holds for some $x_0\in\rd$ and any $t\ge 0$. Then,
equivalently to (\ref{eq:12}), Theorem \ref{thmmain1}
implies Lipschitz contractivity for the transition semigroup
$$(p_tg)(x)\ =\ \int g(z)\, p_t(x,dz)$$
w.r.t.\ the metric $d_f$, i.e.,
\begin{equation}
\label{eqLIP}
\| p_tg\|_{{\rm Lip}(f)}\ \le \ \exp (-ct)\, \| g\|_{{\rm Lip}(f)}
\end{equation}
holds for any $t\ge 0$ and any Lipschitz continuous function $g:\rd\to\r $, where
$$\| g\|_{{\rm Lip}(f)}\ =\ \sup\left\{
\frac{|g(x)-g(y)|}{d_f(x,y)}\, :\, x,y\in\rd\mbox{ s.t.\ }x\neq y\right\} $$
denotes the Lipschitz semi-norm w.r.t.\ $d_f$. An immediate
consequence is the existence of a unique stationary
distribution $\mu $ with finite second moments:

\begin{coroll}[{\bf Convergence to equilibrium}]\label{corSD}
There exists a unique stationary distribution $\mu$ of
$(p_t)_{t\ge 0}$ satisfying
$\int |y|\, \mu (dy) < \infty$ and
\begin{equation}
{\rm Var}_\mu (g) \ \leq \ (2c)^{-1 }\| g\|^2_{{\rm Lip}(f)}\ \mbox{ for any Lipschitz continuous }g:\rd\to\r .\label{eqSD2}
\end{equation}
Moreover, for any probability measure $\nu $ on $\rd$,
\begin{equation}
\label{eqSD3}
W_f(\mu ,\nu p_t)\ \le\ \exp (-ct)\, W_f(\mu ,\nu )\qquad
\mbox{for any }t\ge 0.
\end{equation}
\end{coroll}\medskip

We refer to \cite{BGG,CG} for other recent results on convergence
to equilibrium of diffusion processes in Wasserstein distances.
\smallskip

Further important consequences of (\ref{eqLIP}) are
quantitative non-asymptotic bounds for the decay of
correlations and the bias and variance of ergodic averages.
Let $x_0\in\rd $ and suppose that $(X,\mathbb P)$ is a
solution of (\ref{eq:1}) with initial condition $X_0=x_0$.

\begin{coroll}[{\bf Decay of correlations}]\label{corDC}
For any Lipschitz continuous functions $g,h:\rd\to\r$ and
$s,t\ge 0$,
\begin{equation}\nonumber
{\rm Cov}\,  (g(X_t),h(X_{t+s})) \ \leq \ \frac{1-e^{-2ct}}{2c}\, e^{-cs}\,\| g\|_{{\rm Lip}(f)}\, \| h\|_{{\rm Lip}(f)}.
\end{equation}
\end{coroll}\smallskip

\begin{coroll}[{\bf Bias and variance of ergodic averages}]\label{corEA}
For any Lipschitz continuous function $g:\rd\to\r$ and
$t\in (0,\infty )$,
\begin{eqnarray*}
\left| \mathbb E\left( \frac 1t\int_0^tg(X_s)\,ds\,-\,\int g\,d\mu \right)\right| &\le & \frac{1-e^{-ct}}{ct}\, \| g\|_{{\rm Lip}(f)}\,\int d_f(x_0,y)\,\mu (dy),\ \mbox{ and}\\
{\rm Var}\left( \frac 1t\int_0^tg(X_s)\,ds \right) &\le & \frac{1}{c^2t}\, \| g\|^2_{{\rm Lip}(f)} .
\end{eqnarray*}
\end{coroll}\smallskip

In the variance estimate in Corollary \ref{corEA}, one of
the factors $1/c$ is due to the variance bound (\ref{eqSD2}) w.r.t.\ the stationary distribution, whereas
the second factor $1/c$ bounds the decay rate for the
correlations. Short proofs of Corollaries \ref{corSD}, \ref{corDC},
and \ref{corEA} are included in Section
\ref{secproofs1}. \smallskip

\begin{remark}[{\bf CLT, Gaussian deviation inequality}]\label{remCLT}
The contractivity w.r.t.\ $W_f$ can also be used to prove
a central limit theorem for the ergodic averages
\cite{KW} and a Gaussian deviation inequality strengthening
Corollary \ref{corEA}, cf.\ Remark 2.10 in \cite{Joulin}.
\end{remark}\medskip

\subsection{Examples}\label{sub:firstexamples}

In order to illustrate the quality of the bounds
given in Theorem \ref{thmmain1} and in Corollary \ref{cor16}, we estimate the constant $c$ defined by
(\ref{eq:10}) in different scenarios, and we study the
behaviour of $c$ under perturbations of the
drift $b$.\smallskip

We first consider the situation where $\kappa$ is bounded
from below by a negative constant for any $r$, and by a
positive constant for large $r$:
\begin{lemma}[{\bf Contractivity under lower bounds
on $\protect\boldmath\kappa$}]\label{lemP}
Suppose that
\begin{equation}\label{eq:lowerboundk}
\kappa (r)\geq -L\mbox{ \,for }r\leq R,
\mbox{ \ and \ }\kappa(r)\ge K\mbox{ \,for }r>R
\end{equation} hold with constants $R,L\in[0,\infty)$ and $K\in(0,\infty)$. If $LR_0^2\le 8$ then
\begin{equation}
\alpha^{-1}c^{-1}\ \leq\
\frac{e-1}{2}R^2\,+\,e\sqrt{8K^{-1}}\,R\,+\, 4K^{-1}\ \le\ \frac{3e}{2}\,\max (R^2,8K^{-1}),\label{eq18a}
\end{equation}
and if $LR_0^2\ge 8$ then
\begin{equation}
\alpha^{-1}c^{-1}\ \leq\ 8\sqrt{2\pi}R^{-1}L^{-1/2}(L^{-1}+K^{-1})\exp\left(\frac{LR^2}{8}\right)+32R^{-2}K^{-2}.\label{eq18b}
\end{equation}
\end{lemma}

For diffusions with reflection on a smooth convex domain 
corresponding bounds with $K=\infty$ hold if $R$ is
the diameter of the domain, cf.\ Remark \ref{rem:reflecteddiffusion} above.

\begin{remark}\label{remP}
If $L=0$ then the bound in (\ref{eq18a}) improves to
\begin{equation}
\alpha^{-1}c^{-1}\
 \le\ {2}\,\max (R^2,2K^{-1}).\label{eq18c}
\end{equation}
\end{remark}

The proofs of Lemma \ref{lemP} and Remark \ref{remP} are given in Section \ref{secexamples} below.\medskip

In the first case considered in the lemma, the constant $c$ is at least of order $\min (R^{-2}, K)$. Even if $L=0$ (convex case), this order can not
be improved as one-dimensional Langevin diffusions with potential $U(x)=Kx^2/2$, or, respectively,
with vanishing drift on $(-R/2,R/2)$ demonstrate.
In particular, for $U(x)=Kx^2/2$ with $K>0$, the distance
$W_f$ is equivalent to $W^1$, and the exact decay rate
is $K/2$. This differs from the bounds in (\ref{eq18c})
and (\ref{eq18a}) only by a factor $2$, $6e$ respectively.
% XXXcheck
Thus, {\em if $LR_0^2$ is not too large, the contractivity properties are not affected substantially by non-convexity~!}
\smallskip

In the second case ($LR_0^2\ge 8$), if $K\ge
 \mbox{const.}\cdot L$ then the upper bound for $c^{-1}$ is of order $L^{-3/2}R^{-1}\exp(LR^2/8)$. By the next example, this order in $R$ and $L$ is again optimal:\smallskip

\begin{example}[{\bf Double-well potential with $\protect\boldmath U^{\prime\prime}(x)=-L\mbox{ for }|x|\le R/2$}]\label{exdoublewell}
Consider a Langevin diffusion in $\mathbb R^1$ with a symmetric potential $U\in C^2(\mathbb R)$ satisfying $U(x)=-Lx^2/2$ for $x\in[-R/2,R/2]$, $U^{\prime\prime}\ge -L$, and $\liminf_{|x|\to\infty}U^{\prime\prime}(x)>0$. If $\|\:\cdot\:\|$ is the Euclidean norm then $\kappa(r)= -L$ for $r\in(0,R]$. On the other hand,
let $\tau_0=\inf\{ t\ge 0:X_t=0\}$ denote the first
hitting time of $0$. Then for any initial condition
$x_0>0$,
\begin{equation}
\lim_{t\to\infty}t^{-1}\log\:P_{x_0}[\tau_0>t]\ =\ -\lambda_1(0,\infty)\label{eq:stara}
\end{equation}
where $-\lambda_1(0,\infty)$ is the first Dirichlet
eigenvalue of the generator
$\mathcal Lv=(v''-U'v')/2$ on $(0,\infty)$, cf.~\cite{F}
or see Section \ref{secexamples} below for a short proof
of the corresponding lower bound that is relevant here. If $LR^2\ge 4$ then by
inserting the function $g(x)=\min(\sqrt
Lx,1)$ into the variational characterization of the Dirichlet
eigenvalue, we obtain the upper bound
\begin{equation}
\lambda_1(0,\infty)\
 \le \ \frac 34 e^{1/2} L^{3/2}R\exp (-LR^2/8),\label{eq:star}
\end{equation}
cf.\ Section \ref{secexamples} below.
The estimates (\ref{eq:stara}) and (\ref{eq:star}) seem to indicate that
for $x_0>0$, the Kantorovich distance $W^1(\delta_{-x_0}p_t, \delta_{x_0}p_t)$ decays at most with a rate of
order $L^{3/2}R\exp (-LR^2/8)$. Indeed, under appropriate
growth assumptions on $U(x)$ for $|x|\ge R$, one can prove
that
$$\mathbb P_R\left[\tau_0 >t\right]\ \ge\ 3/4\qquad
\mbox{for any }t\le\lambda_1(0,\infty )^{-1}/4 ,$$
cf.\ Section \ref{secexamples}. Hence for $t\le
3^{-1}e^{-1/2}L^{-3/2}R^{-1}\exp (LR^2/8)$, the
Kantorovich distance $W^1(\delta_Rp_t,\mu )$ between
$\delta_Rp_t$ and the stationary distribution $\mu$ is
bounded from below by a strictly positive constant
that does not depend on $L$ and $R$ if $LR^2\ge 4$.
\end{example}
\bigskip

For analyzing the behaviour of $c$ under perturbations of the drift, we assume that $\| z\| =|\sigma^{-1}z|$ is the
intrinsic metric corresponding to the diffusion matrix,
i.e., $G=(\sigma\sigma^T)^{-1}$. Suppose that
\begin{equation}
\label{eqPA}
b(x)\ =\ b_0(x)+\gamma (x)\qquad\mbox{for any }x\in\r
\end{equation}
with locally Lipschitz continuous functions $b_0,\gamma :
\mathbb{R}^d\to\rd$. For $r>0$ let
\begin{equation}
\label{eqPB}
\kappa_0(r)\ =\ \inf\left\{ -2\,\frac{(x-y)\cdot G(b_0(x)-b_0(y))}{\|x-y\|^2}\, :\,
x,y\in\mathbb R^d\mbox{ s.t.\ }\|x-y\|=r\right\}
\end{equation}
be defined analogously to $\kappa (r)$ with $b$ replaced
by $b_0$. We assume that $\kappa_0$ satisfies the
assumptions (\ref{eqassk}) imposed on $\kappa$ above,
and we define $R_0$ and $R_1$ similarly to (\ref{eqdefr0})
and (\ref{eqdefr1}) but with $\kappa$ replaced by $\kappa_0$.
Now suppose that there exists a constant $R\le R_0$
such that
\begin{equation}
\label{eqPC}
(x-y)\cdot (\gamma (x)-\gamma (y))\ \le\ 0
\qquad\mbox{for any }x,y\in\rd\mbox{ s.t.\ }\| x-y\|\ge R.
\end{equation}
Then $\kappa (r)\ge \kappa_0(r)$ for $r\ge R$, and hence
the constants $R_0$ and $R_1$ defined w.r.t.\ $b$ are
smaller than the corresponding constants defined w.r.t.\
$b_0$. In this situation, we can compare the lower bounds
$c$ and $c_0$ for the contraction rates w.r.t.\ $b$ and
$b_0$ given by (\ref{eq:10}):

\begin{lemma}[{\bf Bounded and Lipschitz perturbations}]\label{lemP1}
Suppose that the drift $b:\rd\to\rd$ is given by (\ref{eqPA}) with $b_0$
and $\gamma$ satisfying the assumptions stated above, and let $c$ and $c_0$ denote the lower bounds for the contraction rates w.r.t.\ $b$ and $b_0$ given by
(\ref{eq:10}).
\begin{enumerate}
\item If $\gamma$ is bounded and (\ref{eqPC}) holds for a constant $R\in [0,R_0]$ then
\begin{equation}
\label{eqPF}
c\ \ge \ c_0\,\exp (-R\sup\|\gamma\| ).
\end{equation}
\item If $\gamma$ satisfies the one-sided Lipschitz
condition
\begin{equation}
\label{eqPG}
(x-y)\cdot G(\gamma (x)-\gamma (y))\ \le\ L\cdot\| x-y\|^2
\qquad\forall\ x,y\in\rd
\end{equation}
with a finite constant $L\in [0,\infty )$
and (\ref{eqPC}) holds for a constant $R\in [0,R_0]$ then
\begin{equation}
\label{eqPH}
c\ \ge \ c_0\,\exp (-LR^2/4 ).
\end{equation}
\end{enumerate}
\end{lemma}

\begin{remark}\label{remP2}
The condition $R\le R_0$ is required in Lemma \ref{lemP1}.
If (\ref{eqPC}) does not hold for $x,y\in\rd $ with
$\| x-y\| \ge R_0$ then the constants $R_0(b)$ and $R_1(b)$
defined w.r.t.\ $b$ are in general greater than the
corresponding constants defined w.r.t.\ $b_0$, i.e., the
region of non-convexity increases by adding the drift $\gamma$. This will also affect the bound in
(\ref{eq:10}) significantly.
\end{remark}

The proof of Lemma \ref{lemP1} is given in Section \ref{secexamples}.

\subsection{Local contractivity and a high-dimensional example}\label{sub:highdim}

Consider again the setup in Section \ref{subsectionRC}.
In some applications, the condition $\liminf_{r\to\infty}
\kappa (r)>0$ imposed above is not satisfied, but the
diffusion process will stay inside a ball $B\subset\rd$ for
a long time with high probability. In this case, one can still prove exponential contractivity up to an error term
that is determined by the exit probabilities from the ball.
Corresponding estimates are useful to prove non-asymptotic
error bounds, i.e., for fixed $t\in (0,\infty )$,
cf.\ e.g.\ \cite{BR, BRHairer, EberleAOP}.\smallskip

Fix $R\in (0,\infty )$ and let $W_{f_R}$ denote the
Kantorovich distance based on the distance function
$d_{f_R}(x,y)=f_R(\| x-y\| )$ given by
\begin{equation}
\label{eqLA}
f_R(r)\ =\ \int_0^r\varphi (s)g_R(s)\, ds\qquad
\mbox{for }r\ge 0,
\end{equation}
where $\varphi$ and $\Phi$ are defined by (\ref{eq:8}), and
\begin{equation}
\label{eqLB}
g_R(r)\ =\ 1-\int_0^{r\wedge R}\frac{\Phi (s)}{\varphi (s)}\, ds\left/ \int_0^{ R}\frac{\Phi (s)}{\varphi (s)}\, ds\right. .
\end{equation}
Notice that
$$g_R(r)=0\qquad\mbox{and}\qquad f_R(r)=f_R(R)
\qquad\mbox{for any }r\ge R,$$
i.e., we have cut the distance at $f_R(R)$.

\begin{thm}[{\bf Local exponential contractivity}]\label{thmLC}
Suppose that the assumptions from Section \ref{subsectionRC} are satisfied except for the condition
$\liminf_{r\to\infty}
\kappa (r)>0$. Then for any $t,R\ge 0$ and
any probability measures $\mu,\nu$ on $\mathbb R^d$,
\begin{eqnarray}\nonumber
W_{f_R}(\mu p_t,\nu p_t) &\leq & \exp ({-c_Rt})\, W_{f_R}(\mu,\nu)\\&& +\, R\cdot \left(\mathbb P_\mu [\tau_{R/2}\le t]
+\mathbb P_\nu [\tau_{R/2}\le t]\right) ,\label{eqLD}
%\ \mbox{\quad and }
%\\
%\nonumber
%W^1(\mu p_t,\nu p_t) & \leq & 2\varphi(R)^{-1}\exp ({-c_Rt})\, W^1(\mu,\nu)\\&&+\,\varphi (R)^{-1}R
%\cdot \left(\mathbb P_\mu [\tau_{R/2}\le t]
%+\mathbb P_\nu [\tau_{R/2}\le t]\right) ,\label{eqLE}
\end{eqnarray}
where $(X_t,\mathbb P_\mu )$ is a diffusion process
satisfying (\ref{eq:1}) with initial distribution $\mu$,
$\tau_{R/2}=\inf\{ t\ge 0:\| X_t\| >R/2\}$ denotes the first
exit time from the ball of radius $R/2$ around $0$, and
\begin{equation}
\label{eqLF}
\frac 1{c_R}\ =\ \alpha\int\limits_0^{R}\Phi(s)\varphi(s)^{-1}\:ds\ =\
\alpha\int\limits_0^{R}\int\limits_0^s\exp\left(\frac 14\int\limits_t^su\kappa(u)^-\:du\right)\:dt\:ds.
\end{equation}
\end{thm}

The proof of the theorem is given in Section \ref{secexamples}. In applications, the exit probabilities
are typically estimated by using appropriate Lyapunov
functions.

\begin{example}[{\bf Stochastic heat equation}]\label{exLE}
We consider the diffusion in $\mathbb R^{d-1}$ given by
$X_t^0\equiv X_t^d\equiv 0$ and
\begin{equation}
\label{eqHE}
dX_t^i\ =\ \left[ d^2\, (X_t^{i+1}-2X_t^i+X_t^{i-1})+
V'(X_t^i)\right]\, dt\, +\,\sqrt d\, dB_t^i,
\end{equation}
$i=1,\ldots ,d-1$, where $V:\r\to\r $ is a $C^2$ function
such that $V''\ge -L$ for a finite constant $L\in\r $.
The equation (\ref{eqHE}) is a spatial discretization
at the grid points $i/d$ ($i=0,1,\ldots ,d$) of the
stochastic heat equation with space-time white noise
and Dirichlet boundary conditions on the interval
$[0,1]$ given by
\begin{equation}
\label{eqHG}
du\ =\ \left(\Delta_{\mbox{Dir}}u\, +\, V'(u)\right)\, dt\, +\, dW
\end{equation}
with the Dirichlet Laplacian $\Delta_{\mbox{Dir}}$ on the
interval $[0,1]$ and a cylindrical Wiener process
$(W_t)_{t\ge 0}$ over the Hilbert space $L^2(0,1)$.
We observe that (\ref{eqHE}) is of the form (\ref{eq:1})
with $\sigma =\sqrt d I_{d-1}$ and $b=-d\nabla U$ where
$$U(x)\ =\ \frac d2 \sum_{i=1}^d\left| x^i-x^{i-1}\right|^2
\, +\, \frac 1d\sum_{i=0}^d V(x^i)$$
for $x=(x^1,\ldots ,x^{d-1})\in\mathbb R^{d-1}$ and
$x^0=x^d=0$. By the discrete Poincar\'e inequality,
$$\sum_{i=1}^d\left| x^i-x^{i-1}\right|^2
\ \ge\ 2\,  (1-\cos (\pi /d))\, \sum_{i=1}^{d-1}\left| x^i\right|^2 .$$
Hence for any $x,\xi\in\mathbb{R}^{d-1}$ and $x^0=x^d
=\xi^0=\xi^d=0$, the lower bound
$$\partial^2_{\xi\xi }U(x)\ =\ d\sum_{i=1}^d\left| \xi^i-\xi^{i-1}\right|^2\, +\,\frac 1d\sum_{i=1}^{d-1}V''(x^i)\left| \xi^i\right|^2
\ \ge\ \frac 1d K_d\sum_{i=1}^{d-1}\left| \xi^i\right|^2$$
holds with
$K_d\ =\ 2\, d^2\, (1-\cos (\pi /d))-L $, and thus
$$(x-y)\cdot (b(x)-b(y))\ =\ -d\, (x-y)\cdot(\nabla U(x)
-\nabla U(y))\ \le\ -K_d\, |x-y|^2$$
for any $x,y\in\mathbb{R}^{d-1}$ where $|\cdot |$ denotes
the Euclidean norm. Choosing for $\|\cdot\|$ the intrinsic
metric $\| x\| =d^{-1/2}|x|$, we obtain
$$\kappa (r)\ \ge\ 2\, K_d\qquad\mbox{ for any }r>0.$$
In particular, the function $\kappa$ is bounded from below uniformly by
a real constant that does not depend on the dimension $d$ since
\begin{equation}
\label{eqHF}
\lim_{d\to\infty} K_d\ =\ \pi^2-L\ >\ -\infty .
\end{equation}
Theorem \ref{thmLC} now shows that for any $R>0$,
local exponential contractivity in the sense of
(\ref{eqLD}) holds on the ball
$$B_{R/2}\ =\ \{ x\in\mathbb{R}^{d-1}: \| x\|\le R/2\}
\ =\ \{ x\in\mathbb{R}^{d-1}: | x|\le d^{1/2}R/2\} $$
with rate $c_R$ satisfying
\begin{eqnarray*}
\frac{1}{c_R}&\le &4\sqrt{\pi }R^{-1}|K_d|^{-3/2}\exp (-K_dR^2/4)\qquad\mbox{for }K_dR^2\le -4,\\
\frac{1}{c_R}&\le & (e-1)R^2/2\qquad\qquad \mbox{for }-4\le K_dR^2<0,\\
\frac{1}{c_R}&\le & R^2/2\qquad\qquad \mbox{for }K_d=0
\mbox{ respectively.}
\end{eqnarray*}
Here the explicit upper bounds are obtained analogously 
as in the proof of Lemma \ref{lemP}.
For $K_d>0$, strict convexity holds, and we obtain global
exponential contractivity with a dimension-independent
rate. We remark that because of (\ref{eqHF}), the bounds
also carry over to the limiting SPDE (\ref{eqHG}) for which
they imply local exponential contractivity on balls w.r.t.\ the $L^2$
norm.
\end{example}

\section{Main results for componentwise reflection couplings}\label{secresultsprod}

\subsection{Componentwise reflection couplings and contractivity on product spaces}\label{subsectioncrc}

We now consider a system
\begin{equation}
\label{eqcircstar}
dX_t^i\ =\ b^i(X_t)\, dt\, +\, dB_t^i ,\qquad
i=1,\ldots ,n,
\end{equation}
of $n$ interacting diffusion processes taking values
in $\mathbb{R}^{d_i}$, $d_i\in\mathbb{N}$. Here
$B^i$, $i=1,\ldots ,n$, are independent Brownian motions in
$\mathbb{R}^{d_i}$, $X=(X^1,\ldots ,X^n)$ is a diffusion
process taking values in $\rd$ where $d=\sum_{i=1}^nd_i$, and $b^i:\rd\to
\mathbb{R}^{d_i}$ are locally Lipschitz continuous
functions. We will assume that
\begin{equation}
\label{eq14star}
b^i(x)\ =\ b_0^i(x^i)\, +\, \gamma^i(x) ,\qquad
i=1,\ldots ,n,
\end{equation}
where the functions $b_0^i:\mathbb{R}^{d_i}\to\mathbb{R}^{d_i}$ are
locally Lipschitz continuous, and $\gamma^i:\rd\to
\mathbb{R}^{d_i}$ are ``sufficiently small'' perturbations,
cf.\ Theorem \ref{thmA} below. In particular, for $\gamma^i
\equiv 0$ the components $X^1,\ldots ,X^n$ are independent.
\smallskip

To analyse contraction properties of the process $X$, one could use a reflection coupling on $\rd$ and apply the results above based on a distance function of the form
$d_f(x,y)=f(|x-y|)$. In some applications, this approach
does indeed provide dimension-free bounds, cf.\ Example
\ref{exLE} above. However, in the product case $\gamma^i
\equiv 0$ it leads in general to lower bounds for contraction rates that degenerate rapidly as $n\to\infty$,
even though one would expect exponential contractivity
with the minimum of the contraction rates for the components. The reason is that the approach requires
convexity outside a Euclidean ball in $\rd$ whereas
in corresponding product models, in general convexity only holds if all components are outside given balls in
$\mathbb{R}^{d_i}$.\smallskip

Instead, we now consider contractivity w.r.t.\ Kantorovich
distances $W_{f,w}$ based on distance functions on
$\rd=\mathbb{R}^{d_1+\cdots +d_n}$ of the form
\begin{equation}
\label{eq14delta}
d_{f,w}(x,y)\ =\ \sum_{i=1}^nf_i(|x^i-y^i|)\, w_i\, .
\end{equation}
Here $f_i:[0,\infty )\to [0,\infty )$, $1\le i\le n$,
are strictly increasing concave $C^1$ functions
with $f_i(0)=0$ and $f_i^\prime (0)=1$ that are obtained
from $b_0^i$ in the same way as $f$ has been obtained
from $b$ above, and $w_i\in (0,1]$ are positive weights.
In many applications, one can choose $w_i=1$ for any $i$.
The corresponding distance will then be denoted by $d_{1,f}$. Notice that $d_{1,f}$ is bounded from above
by the $\ell^1$ distance
$$d_{\ell^1}(x,y)\ =\ \sum_{i=1}^n|x^i-y^i| .$$
Hence $W_{1,f}$ is bounded from above by the Kantorovich
distance $W_{\ell^1}$ based on $d_{\ell^1}$.\medskip

For $r\in (0,\infty )$ let
\begin{equation}\label{eqKI}
\kappa_i(r)\:=\:r^{-2}\, \inf\left\{ -2\,(x-y)\cdot (b_0^i(x)-b_0^i(y))\, :\,
x,y\in\mathbb R^d\mbox{ s.t.\ } |x-y|=r\right\} .
\end{equation}
Similarly as above, we assume that for $1\le i\le n$,
\begin{equation}\label{eqAK1}
\kappa_i:(0, \infty )\to\r \mbox{ is continuous with }\,
\liminf_{r\to\infty}\kappa_i(r)>0.
\end{equation}
Moreover, we assume
\begin{equation}\label{eqAK2}
\lim_{r\to 0}r\kappa_i(r)=0.
\end{equation}
Let $R_0^i$, $R_1^i$, $g_i(r)$, $\varphi_i(r)$,
$f_i(r)$ and $\Phi_i(r)=\int_0^r\varphi_i(s)\, ds$
be defined analogously to (\ref{eqdefr0}), (\ref{eqdefr1})
and (\ref{eq:8}) with $\kappa$ replaced by $\kappa_i$.
Moreover, we define $c_i\in (0,\infty )$ by
\begin{equation}
\label{eq14dot}
\frac 1{c_i}\ =\ \int\limits_0^{R_1^i}\Phi_i(s)\varphi_i(s)^{-1}\:ds\ =\
\int\limits_0^{R_1^i}\int\limits_0^s\exp\left(\frac 14\int\limits_t^su\kappa_i(u)^-\:du\right)\:dt\:ds\, .
\end{equation}
Recall that by Theorem \ref{thmmain1} and Corollary \ref{cor16},
$c_i$ is a lower bound for the contraction rate of the
diffusion process $\tilde X^i$ on $\mathbb{R}^{d_i}$
satisfying the s.d.e.\ $d\tilde X_t^i=b_0^i(\tilde X_t^i)\,dt\,+\,dB_t^i$.
\smallskip

 Let $p_t(x,dy)$ denote the transition kernels of the
 diffusion process $X_t=(X_t^1,\ldots ,X_t^d)$ on $\rd$
 satisfying (\ref{eqcircstar}). We now state our second
 main result:

\begin{thm}[{\bf Exponential contractivity on
product spaces}]\label{thmA}
Suppose that (\ref{eqAK1}) and (\ref{eqAK2}) hold, and
suppose that there exist constants $\varepsilon_i\in [0,c_i)$, $1\le i\le n$, such
that for any $x,y\in\rd $,
\begin{equation}
\label{eq14dotdot}
\sum_{i=1}^n|\gamma^i(x)-\gamma^i(y)|\, w_i\ \le\
\sum_{i=1 }^n\varepsilon_i\, f_i(|x^i-y^i|)\, w_i .
\end{equation}
Then for any $t\geq 0$ and any probability measures $\mu,\nu$ on $\mathbb R^d$,
\begin{eqnarray}
W_{f,w}(\mu p_t,\nu p_t) &\leq & \exp ({-ct})\, W_{f,w}(\mu,\nu),\ \mbox{\qquad and }\label{eqAX}\\
W_{\ell^1}(\mu p_t,\nu p_t) & \leq & A\, \exp ({-ct})\, W_{\ell^1}(\mu,\nu),\label{eqAY}
\end{eqnarray}
where $\ c\ =\ \min\limits_{i=1,\ldots ,n} (c_i-\varepsilon_i)\qquad
\mbox{and}\qquad A\ =\ 2\left/\min\limits_{i=1,\ldots ,n}
(\varphi_i(R_0^i)w_i)\right. .$
\end{thm}

\begin{example}[{\bf Product model}]\label{exProduct}
In the product case, $\gamma^i\equiv 0$ for any $i$.
Hence Condition (\ref{eq14dotdot}) is satisfied with
$\varepsilon_i=0$, and, therefore,
$$ W_{f,w}(\mu p_t,\nu p_t) \ \leq \ \exp ({-ct})\, W_{f,w}(\mu,\nu)$$
holds with $c=\min c_i$ for any choice of the weights
$w_1,\ldots , w_n$.
\end{example}

More generally than in the example, suppose now that
$\gamma =(\gamma^1,\ldots ,\gamma^n)$ satisfies an
$\ell^1$-Lipschitz condition
\begin{equation}
\label{eqA}
\sum_{i=1}^n|\gamma^i(x)-\gamma^i(y)|\ \le\
\lambda\, \sum_{i=1 }^n |x^i-y^i|\qquad
\forall\ x,y\in\rd .
\end{equation}
Then exponential contractivity holds for the perturbed
product model provided $\lambda <c_i\varphi (R_0^i)/2$
for any $i$:

\begin{coroll}[{\bf Perturbations of product models}]\label{corB}
Suppose that (\ref{eq14star}), (\ref{eqAK1}), (\ref{eqAK2}) and (\ref{eqA}) hold with
$\lambda\in [0,\infty )$. Then
for any $t\geq 0$ and any probability measures $\mu,\nu$ on $\mathbb R^d$,
\begin{eqnarray}
\label{eqB}
W_{f,1}(\mu p_t,\nu p_t) & \leq &\exp ({-ct})\, W_{f,1}(\mu,\nu),\ \mbox{ and }\\
W_{\ell^1}(\mu p_t,\nu p_t) & \leq &  A\exp ({-ct})\, W_{\ell^1}(\mu,\nu),
\end{eqnarray}
where $\ c=\min\limits_{i=1,\ldots n}(c_i-2\lambda
\varphi_i(R_0^i)^{-1})\ $ and $\ A=2\max\limits_{i=1,\ldots n}\varphi_i(R_0^i)^{-1}$.
\end{coroll}

The inituitive idea of proof for Theorem \ref{thmA} is to
construct a coupling $(X_t,Y_t)$ of two solutions of
(\ref{eqcircstar}) by applying a reflection coupling
individually for each component $(X_t^i,Y_t^i)$ if
$X_t^i\neq Y_t^i$, and a synchronuous coupling if $X_t^i=Y_t^i$.
In the product case this just means that $X_t^i=Y_t^i$
for any $t\ge \tau^i$ where $\tau^i=\inf\{
t\ge 0: X_t^i=Y_t^i\} $ is the coupling time for the
$i$-th component. In the non-product case, however,
$X_t^i$ and $Y_t^i$ can move apart again after the time
$\tau^i$ due to interactions with other components. In that
case it is not clear how to define a coupling as described
above rigorously. Instead we will use a regularized version
where reflection coupling is applied to the $i$-th
component whenever $|X_t^i-Y_t^i|\ge\delta $ for a given
constant $\delta >0$, and synchronuous coupling is applied
whenever $|X_t^i-Y_t^i|\le\delta /2 $. A precise
description of the coupling and the proofs of Theorem
\ref{thmA} and Corollary \ref{corB} are given in Sections
\ref{secprod} and \ref{secinteract} below.

\subsection{Consequences}

The contractivity results in Theorem \ref{thmA} and
Corollary \ref{corB} have
corresponding consequences as the contractivity results
in the non-product case, cf.\ Section \ref{sub:cons}
above. An important difference to be noted is, however,
that on product spaces,
$$d_{f,w}(x,y)\ \le\ \sum_{i=1}^n|x^i-y^i|\ \le\
n^{1/2}\, |x-y|$$
by the Cauchy-Schwarz inequality. Therefore, an
additional factor $n$ occurs in the variance bounds
from Corollaries \ref{corSD}, \ref{corDC} and
\ref{corEA} on product spaces. Apart from this
additional factor, all results in Section \ref{sub:cons}
carry over to the setup considered in Section
\ref{subsectioncrc}.

\subsection{Interacting Langevin diffusions}

As an illustration of the results in Section \ref{subsectioncrc}, we consider a system
\begin{equation}
\label{eqE}
dX_t^i\ =\ -\frac 12\nabla U(X_t^i)\, dt\,-\,
\sum_{j=1}^na_{ij}\, \nabla V(X_t^i-X_t^j)\, dt \, +\, dB_t^i
\end{equation}
of $n$ interacting overdamped Langevin diffusions taking
values in $\mathbb{R}^k$ for some $k\in\mathbb{N}$.
Here $B^1,\ldots ,B^n$ are independent Brownian motions
in $\mathbb{R}^k$, $U\in C^2(\mathbb{R}^k)$ is strictly
convex outside a given ball, the interaction potential
$V$ is in $C^2(\mathbb{R}^k)$ with bounded second
derivatives, and $a_{ij}$, $1\le i,j\le n$, are
finite real constants. For example, we are interested in
{\em nearest-neighbour interactions} and {\em mean-field
interactions} given by
\begin{eqnarray}
\label{eqC}
a_{ij}& =& \left\{
\begin{array}{ll}\alpha /2\ \ &\mbox{if }i-j\equiv 1
\mbox{ mod }n\mbox{ \ or \ }i-j\equiv -1
\mbox{ mod }n,\\
0&\mbox{otherwise,}
\end{array}\right. \\
\label{eqD}
a_{ij}& =&\alpha\, n^{-1}\qquad\mbox{respectively,}
\end{eqnarray}
where $\alpha\in\r $ is a finite coupling constant.
\smallskip

Choosing $b_0^i(x^i)=-\nabla U(x^i)/2$ and $\gamma^i(x)=
-\sum_{j=1}^na_{ij}\nabla V(x^i-x^j)$, we observe that the
function
$$\kappa_i(r)\ =\ \inf\left\{\int\nolimits_0^1\partial^2_{(x-y)/|x-y|}U((1-t)x+ty)\,
dt :x,y\in\mathbb R^k \mbox{ s.t.}\, |x-y| =r\right\}$$
does not depend on $i$. Let $\varphi$ and $f$ be the
corresponding functions given by (\ref{eq:8}), and consider
the distance
$$d_{1,f}(x,y)\ =\ \sum_{i=1}^nf(|x^i-y^i|).$$
Morover, let $c$ be given by (\ref{eq:10}) with $\alpha =1$, i.e., $c$ is the lower bound for the contraction rate
of the diffusion process $Y$ in $\mathbb R^k$ satisfying
$dY=-\frac 12\nabla U(Y)\, dt\, +\, dB$. We note that
$\gamma $ satisfies the $\ell^1$ Lipschitz condition
(\ref{eqA}) with
$$\lambda\ =\ M\cdot\max_i\sum_{j=1}^n \left(
|a_{ij}|+|a_{ji}|\right)$$
where $M=\sup\|\nabla^2V\| $. Therefore, if
$$\sum_{j=1}^n\left(
|a_{ij}|+|a_{ji}|\right)\ \le\ c\, \varphi (R_0)\, M^{-1}
$$
then by Corollary \ref{corB}, contractivity in the sense
of (\ref{eqB}) holds with contraction rate
$$\bar c\ =\ c-2\lambda\varphi (R_0)^{-1}\ >\ 0. $$
In particular, in the nearest neighbour and mean field case,
we obtain contractivity with a rate that does not depend
on the dimension if $ \alpha$ is small:

\begin{coroll}[{\bf Mean field and nearest neighbour interactions}]\label{corC}
Let $p_t$, $t\ge 0$, denote the transition kernels
of the diffusion process on $\mathbb R^{nk}$ solving
(\ref{eqE}). Suppose that $\sup\|\nabla^2V\| <\infty $
and that $a_{ij}$ is given by (\ref{eqC}) or by (\ref{eqD})
with $\alpha\in\r $. Then there exist finite constants
$c,\theta ,A\in (0,\infty )$ that do not depend on the
dimension $n$ such that
\begin{eqnarray}
\label{eqF}
W_{f,1}(\mu p_t,\nu p_t) & \leq & e^{(\theta\alpha -c)t}\, W_{f,1}(\mu,\nu),\ \mbox{ and }\\
\label{eqG}
W_{\ell^1}(\mu p_t,\nu p_t) & \leq &  A\, e^{(\theta\alpha -c)t}\, W_{\ell^1}(\mu,\nu),
\end{eqnarray}
hold for any $t\geq 0$ and any probability measures $\mu,\nu$ on $\mathbb R^{nk}$. In particular, exponential contractivity holds for $\alpha <c/\theta$.
\end{coroll}

The bounds in (\ref{eqF}) and (\ref{eqG}) are not sharp.
However, it is known that for example in mean field models
where $U$ is a double-well potential and $V$ is quadratic,
exponential contractivity with a rate independent of the
dimension can not be expected to hold for large $\alpha$.
Indeed, in this case the corresponding McKean-Vlasov
process has several stationary distributions if
$\alpha > \alpha_1$ for some critical parameter
$\alpha_1\in (0,\infty )$, cf.\ \cite{HerrmannTugaut, HerrmannTugauta}.

%\subsection{Outlook}
%%XXX Central Limit Theorem, non-asymptotic variance estimates
%%
%%The contraction properties in Kantorovich distances can be applied to derive non-asymptotic variance bounds
%%for ergodic averages of the form $t^{-1}\int_0^tf(X_s)\, ds$
%%where $X$ is a diffusion process satisfying (\ref{eq:1}).
%%Corresponding results in discrete time have been proven
%%by Joulin and Ollivier in \cite{JOXXX}. In particular,
%%under appropriate assumptions, Wasserstein contractivity
%%implies a central limit theorem for ergodic averages,
%%cf.~Komorowski and Walczuk \cite{KWXXX}.
%
%XXX Generalization to discrete time and continuous state
%space; algorithmic applications
%
%XXX Combination with Lyapunov condition

\section{Proofs for Reflection Coupling}\label{secproofs1}

In this section, we first motivate our particular choice
of the function $f$, and we prove Theorem \ref{thmmain1}.
Afterwards, we prove Corollaries \ref{cor16},
\ref{corSD}, \ref{corDC} and \ref{corEA}.\medskip

Let $r_t=\|X_t-Y_t\|$ where $(X,Y)$ is a reflection coupling of two solutions of (\ref{eq:1}). Our goal
is to find an explicit concave increasing function
$f:[0,\infty )\to [0,\infty )$ with $f(0)=0$ and
$f^\prime (0)=1$ such that $e^{ct}f(r_t)$ is a (local)
supermartingale for $t$ less than the coupling time $T$
with a constant $c>0$ that we are trying to maximize by
the choice of $f$.\smallskip

An application of It\^o's formula to the s.d.e.\ (\ref{eq:3}) satisfied by the difference process $Z_t=X_t-Y_t$ shows that the following It\^o equations hold
almost surely for $t<T$ whenever $f$ is $C^1$ and $f^\prime$ is absolutely continuous:
\begin{eqnarray}
\nonumber
d\| Z_t\|^2&=&4\, |\sigma^{-1}Z_t|^{-1}\| Z_t\|^2\:dW_t\\
&&+\, 2\, Z_t\cdot G(b(X_t)-b(Y_t))\:dt\, +\,4\, |\sigma^{-1}Z_t|^{-2}\| Z_t\|^2\:dt,\nonumber\\
\nonumber
dr_t&=&2\, |\sigma^{-1}Z_t|^{-1}r_t\:dW_t\,+
\,r_t^{-1}Z_t\cdot G(b(X_t)-b(Y_t))\:dt,\ \ \mbox{ and}\\
\nonumber
df(r_t)&=&2\, |\sigma^{-1}Z_t|^{-1}r_t\, f'(r_t)\:dW_t\\
\label{eq:13}&&+\, r_t^{-1}Z_t\cdot G(b(X_t)-b(Y_t))f'(r_t)\:dt\, +\, 2\, |\sigma^{-1}Z_t|^{-2}r_t^2f''(r_t)\:dt.
\end{eqnarray}
By definition of the function $\kappa $, the drift term on
the right hand side of (\ref{eq:13}) is bounded from above by
\begin{equation}\label{eq21a}
\beta_t\::=\:2\, |\sigma^{-1}Z_t|^{-2}r_t^2\cdot\left(f''(r_t)-\frac 14\, r_t\, \kappa(r_t)f'(r_t)\right).
\end{equation}
Hence the process $e^{ct}f(r_t)$ is a supermartingale for $t<T$ if $\beta_t\le -cf(r_t)$. Since
\begin{equation}
\label{21b}
|\sigma^{-1}z|^2\ \le\ \alpha\| z\|^2\qquad
\mbox{for any }z\in\mathbb{R}^d
\end{equation}
with $\alpha$ defined as in Theorem \ref{thmmain1}, a sufficient condition is
\begin{equation}
\label{21c}
f''(r)-\frac 14r\kappa (r)f'(r)\ \le\ -\frac{\alpha c}{2}f(r)\qquad
\mbox{for a.e.\ }r>0.
\end{equation}
We now first observe that this equation holds with $c=0$
(i.e., $f(r_t)$ is a supermartingale for $t<T$) if $f$ is chosen such that $f'(r)=\varphi (r)=\exp
(-\int_0^rs\kappa (s)^- ds/4)$. Indeed, $f(r)=\int_0^r
\varphi (s)\, ds$ is the least concave among all concave
functions $f$ satisfying $\beta_t\le 0$.\smallskip

To satisfy the stronger condition $\beta_t\le -cf(r_t)$
with $c>0$, we make the ansatz
\begin{equation}
\label{21f}
f'(r)\ =\ \varphi (r)\, g(r)
\end{equation}
with a decreasing absolutely continuous function $g\ge 1/2$
such that $g(0)=1$. Note that the condition $g\ge 0$ is
required to ensure that $f$ is non-decreasing. By replacing
this condition by the stronger condition $g\ge 1/2$,
we are loosing at most a factor $2$ in the estimates below.
On the other hand, the condition $1/2\le g\le 1$ has the
huge advantage of ensuring that
\begin{equation}
\label{21e}
\Phi /2\ \le\ f\ \le \ \Phi
\end{equation}
where $\Phi (r)=\int_0^r\varphi (s)\, ds$. The ansatz
(\ref{21f}) yields
$$f''\ =\ -\frac 14\, r\kappa^-f\,+\,\varphi g'\
\le\ \frac 14\, r\kappa f\,+\,\varphi g' ,$$
i.e., Condition (\ref{21c}) is satisfied if
\begin{equation}
\label{21d}
g'\ \le\ -\frac{\alpha c}{2}f/\varphi\, .\qquad\mbox{almost surely.}
\end{equation}
We will see in the proof below that for $r\ge R_1$, Condition (\ref{21c}) is automatically satisfied since
$\kappa $ is sufficiently positive. Therefore, it is
enough to assume that (\ref{21d}) holds on $(0,R_1)$.
\smallskip

Now on the one hand, if (\ref{21d}) is satisfied on $(0,R_1)$
then
$$g(R_1)\ \le\ 1-\frac{\alpha c}{2}
\int_0^{R_1}f(s)\varphi (s)^{-1}\, ds\ \le\ 1-\frac{\alpha c}{4}\int_0^{R_1}\Phi (s) \varphi (s)^{-1}\, ds .$$
This condition can only be satisfied with a function $g$
taking values in $[1/2,1]$ if
$$\alpha\, c\ \le\ 2\left/\int_0^{R_1}\Phi (s) \varphi (s)^{-1}\, ds \right. .$$
On the other hand, by choosing
\begin{equation}\label{21g}
g'(r)\ =\ -\frac{\Phi(r)}{2\varphi (r)} \left/\int\nolimits_0^{R_1}\frac{\Phi (s)}{\varphi (s)}\,ds\right.\qquad\mbox{for }r<R_1,
\end{equation}
Condition (\ref{21d}) is satisfied with the constant
$$\alpha\, c\ =\ 1\left/\int_0^{R_1}\Phi (s) \varphi (s)^{-1}\, ds \right. .$$
This shows that up to a factor 2, choosing $g$ as in (\ref{21g}) is the best we can do under the assumptions that we have made.\medskip

The considerations above explain the particular choice of the function $f$ made in (\ref{eq:8}). Once this choice has
been made, the proof of Theorem \ref{thmmain1} is
almost straightforward:\bigskip

{\bf Proof of Theorem \ref{thmmain1}.}
As remarked above, the drift in the s.d.e.\ (\ref{eq:13})
for $f(r_t)$ is bounded from above by $\beta_t$ defined
by (\ref{eq21a}). We now show that by our choice of $f$
in (\ref{eq:8}), this expression is
%. By   and ,
%\begin{equation}
%%\label{eq:13}
%df(r_t)=2|\sigma^{-1}Z_t|^{-1}r_tf'(r_t)\,dW_t+r_t^{-1}Z_t\cdot G(b(X_t)-b(Y_t))f'(r_t)\:dt+2|\sigma^{-1}Z_t|^{-2}r_t^2f''(r_t)\,dt
%\end{equation}
%a.s.\ for $t<T$. The drift is bounded from above by
%$\beta_t:=2|\sigma^{-1}Z_t|^{-2}r_t^2\left(f''(r_t)-r_t\kappa(r_t)f'(r_t)/4\right)$.
%We show that by our choice of $f$, this expression is
smaller than $-cf(r_t)$ where $c$ is given by
(\ref{eq:10}). Indeed, for $r<R_1$,
\begin{eqnarray}
f''(r)&=& -\frac 14 r\kappa(r)^-\varphi(r)g(r)-\frac 12\Phi(r)\left/\int\limits_0^{R_1}\Phi(s)\varphi(s)^{-1}\:ds\right.\label{eq:14}\\
&\leq &\ \frac 14 r\kappa(r)f'(r)-\frac 12 f(r)\left/\int\limits_0^{R_1}\Phi(s)\varphi(s)^{-1}\:ds\right..\nonumber
\end{eqnarray}
%\begin{eqnarray}\label{eq:14}
%f''(r) &=& -\frac 14 r\kappa(r)^-\varphi(r)g(r)-\frac 12\Phi(r)\left/\int\nolimits_0^{R_1}\frac{\Phi (s)}{\varphi (s)}\,ds\right.\\
%\leq\,\frac 14 r\kappa(r)f'(r)-\frac 12 f(r)\left/\int\nolimits_0^{R_1}\frac {\Phi (s)}{\varphi (s)}\,ds\right..
%\end{eqnarray}
For $r > R_1$, we have $f'(r)=\varphi(r)/2=\varphi(R_0)/2$ and
$\kappa(r)R_1(R_1-R_0)\geq 8$ by definition of $R_1$, whence
\begin{eqnarray}\nonumber
f''(r)-\frac 14 r\kappa(r)f'(r)& = &-\frac 18
r\kappa(r)\varphi(R_0)\ \leq\
-\frac{\varphi(R_0)}{R_1-R_0}\cdot\frac{r}{R_1}\\
\label{eq:15}
&\leq & -\frac{\varphi(R_0)}{R_1-R_0}\cdot\frac{\Phi(r)}{\Phi(R_1)}\
\leq\ -\frac
12\Phi(r)\left/\int\nolimits_{R_0}^{R_1}\Phi(s)\varphi(s)^{-1}\:ds\right.\\
\nonumber &\leq & -\frac
12f(r)\left/\int\nolimits_0^{R_1}\Phi(s)\varphi(s)^{-1}\:ds\right..
\end{eqnarray}
Here we have used that for $r\ge R_0$, the function $\varphi (r)$ is constant, and, therefore, $\Phi(r)\ =\ \Phi(R_0)+(r-R_0)\, \varphi(R_0)$,
and
\begin{eqnarray*}
\int\nolimits_{R_0}^{R_1}\Phi(s)\varphi(s)^{-1}\:ds &=&\int\nolimits_{R_0}^{R_1}(\Phi(R_0)+(s-R_0)\varphi(R_0))\varphi(R_0)^{-1}\,ds\\
&=&{\Phi(R_0)}{\varphi(R_0)^{-1}}(R_1-R_0)+(R_1-R_0)^2/2\\
&\geq &(R_1-R_0)\left(\Phi(R_0)+(R_1-R_0)\varphi(R_0)\right)\varphi(R_0)^{-1}/2\\
&= &(R_1-R_0)\Phi(R_1)\varphi(R_0)^{-1}/2.
\end{eqnarray*}
By (\ref{eq:14}) and (\ref{eq:15}), we conclude that $\beta_t\leq -cf(r_t)$. Optional stopping in (\ref{eq:13}) at $T_k=\inf\{t\geq 0:\, r_t\not\in(k^{-1},k)\}$ now implies
$$\E[f(r_t)\:;\:t<T_k]\ \leq\ -c\int\nolimits_0^t\E[f(r_s)\:;\:s<T_k]\:ds\, $$ for any $k\in\mathbb N$ and $t\geq 0$.
The assertion follows for $k\to\infty$ since $r_t=0$ for $t\geq T$, and $T=\sup T_k$ by non-explosiveness.
$\Box$\medskip

{\bf Proof of Corollary \ref{cor16}.} Let $(X,Y)$ be a
reflection coupling of two solutions of (\ref{eq:1})
with joint initial distribution $(X_0,Y_0)\sim\eta $.
Then by Theorem \ref{thmmain1},
\begin{eqnarray*}
W_f(\mu p_t,\nu p_t) & \le & \E \left[ d_f(X_t,Y_t)\right]
\ \le \ e^{-ct}\,  \E \left[ d_f(X_0,Y_0)\right]\\
&=& e^{-ct}\,\int d_f(x,y)\,\eta (dx\, dy)
\end{eqnarray*}
for any $t\ge 0$. The estimate (\ref{eq:12}) now follows
by taking the infimum over all couplings $\eta$ of two
given probability measures $\mu$ and $\nu$ on
$\mathbb R^d$. Moreover, (\ref{eq:12a}) follows from
(\ref{eq:12}) by (\ref{eq113a}).
$\Box$\bigskip

%{\bf Proof of Lemma \ref{lemP}.}
% since $\varphi(r)=\varphi(R_0)$ and $\Phi(r)=\Phi(R_0)+(r-R_0)\varphi(R_0)$ for $r\geq R_0$,
%\begin{equation}
%\label{eq:17}
%\alpha^{-1}c^{-1}\:=\:\int\nolimits_0^{R_1}\Phi(s)\varphi(s)^{-1}\:ds\:=\:\int\nolimits_0^{R_0}\Phi(s)\varphi(s)^{-1}\:ds+(R_1-R_0)\Phi(R_0)\varphi(R_0)^{-1}+(R_1-R_0)^2/2.
%\end{equation}
%The lower bounds on the function $\kappa$ imply the upper bounds $R_0\leq R$, $R_1-R_0\leq\min(8/(KR_0),\sqrt{8/K})$,\\
%$
%\Phi(r)\varphi(r)^{-1}\,\leq\,\int\nolimits_0^r\exp(L(r^2-t^2)/8)\,dt\,
%\leq\, \min(\sqrt{2\pi/L},r)\exp(Lr^2/8)\,\mbox{ for }r\leq R_0,\ \mbox{ and}$
%$$\int\limits_0^{R_0}\Phi(r)\varphi(r)^{-1}\:dr\ \leq\ \left\{ \begin{array}{ll}
%4L^{-1}(\exp (LR_0^2/8)-1)\ \le\ (e-1)R_0^2/2 &\qquad\mbox{ if }LR_0^2/8\leq1,\\
%8\sqrt{2\pi}L^{-3/2}R_0^{-1}\exp(LR_0^2/8)&\qquad\mbox{ if }LR_0^2/8\geq1.
%\end{array}\right.
%$$
%Combining these estimates, we obtain by (\ref{eq:17}),
%$$
%\alpha^{-1}c^{-1}\ \leq\
%\left\{ \begin{array}{ll}
%(e-1)R^2/2+e\sqrt{8/K}R+4/K\ \le\ (3e/2)\max (R^2,8/K)&\qquad\mbox{ if }LR_0^2/8\:\leq\:1,\\
%8\sqrt{2\pi}R^{-1}L^{-1/2}(L^{-1}+K^{-1})\exp(LR^2/8)+32R^{-2}K^{-2}&\qquad\mbox{ if }LR_0^2/8\geq 1.\end{array}\right.\Box
%$$
%\medskip

Next, we are going to prove the results in Section
\ref{sub:cons}. Suppose that (\ref{eqFM}) holds,
$\| z\| =|\sigma^{-1}z|$ is the intrinsic metric,
and $b$ is in $C^1$. Corollary \ref{cor16} implies
$$\int |y|\, p_t(x,dy)\ \le\ \int |y|\, p_t(x_0,dy)\,
+\, W^1(p_t(x,\cdot ),p_t(x_0,\cdot ))\, <\,\infty $$
for any $t\ge 0$ and any $x\in\rd $. In particular,
$(p_tg)(x)=\int g(y)\,p_t(x,dy)$ is defined for any
Lipschitz continuous function $g:\rd\to\r $, and
$$|(p_tg)(x)-(p_tg)(y)|\ =\ \left| \mathbb E [g(X_t)-
g(Y_t)]\right|\ \le\ \| g\|_{{\rm Lip}(f)}
\mathbb E [d_f(X_t,Y_t)]$$
for any coupling $(X_t,Y_t)$ of $p_t(x,\cdot )$ and
$p_t(y,\cdot )$. Hence by Theorem \ref{thmmain1},
\begin{equation}
\label{eqLIPA}
|(p_tg)(x)-(p_tg)(y)|\ \le\  e^{-ct}\,\| g\|_{{\rm Lip}(f)}
\, d_f(x,y),
\end{equation}
i.e., $p_t$ satisfies the exponential contractivity
condition (\ref{eqLIP}) w.r.t.\ $\| \cdot\|_{{\rm Lip}(f)}$. If $p_tg$ is $C^1$ then by (\ref{eqLIPA}) and
since
$$d_f(x,y)\ \le\ \| x-y\| \ =\ |\sigma^{-1}(x-y)|
\qquad\forall\,  x,y\in\rd ,$$
we obtain the uniform gradient bound
\begin{equation}
\label{eqLIPB}
\sup\left|\sigma^T\nabla p_tg\right|\ \le\  e^{-ct}\,\| g\|_{{\rm Lip}(f)}
\qquad\forall\ t\ge 0.
\end{equation}
It is well-known that this bound can be used to control
variances w.r.t.\ the measures $p_t(x,\cdot )$:

\begin{lemma}\label{lemV}
For any $t\ge 0$, $x\in\rd$, and any Lipschitz continuous
$g:\rd\to\r$,
\begin{equation}
\label{eqVA}
{\rm Var}_{p_t(x,\cdot )}(g)\ \le\  \frac{1-\exp ({-2ct})}{2c} \,\| g\|^2_{{\rm Lip}(f)}
.
\end{equation}
\end{lemma}

{\bf Proof.} We may assume $g\in C^2(\rd )$ and $t>0$.
Then, by standard elliptic regularity results,
$(t,x)\mapsto (p_tg)(x)$ is differentiable in $t$ and $x$,
and $$\frac{d}{dt}\, p_tg\ =\ \mathcal Lp_tg\ =\ p_t\mathcal Lg$$
where $\mathcal L=\frac 12\sum a_{ij}\frac{\partial^2}{\partial x^i\partial x^j}+b(x)\cdot\nabla $,
$a=\sigma\sigma^T$, is the generator of $(X_t)$, cf.\
e.g.\ \cite{Stroock, Roy}. In particular, for
$s\in (0,t)$,
\begin{eqnarray*}
\frac{d}{ds}\, p_s(p_{t-s}g)^2 &=& p_s\,
\left(\mathcal L(p_{t-s}g)^2-2p_{t-s}g\, \mathcal Lp_{t-s}g\right)\\
&=& p_s\,\left|\sigma^T\nabla p_{t-s}g\right|^2
\ \le\ e^{-2c(t-s)}\| g\|^2_{{\rm Lip}(f)}
\end{eqnarray*}
by (\ref{eqLIPB}). Integrating w.r.t.\ $s$, we obtain
$$p_tg^2-(p_tg)^2\ \le\ \frac{1-\exp ({-2ct})}{2c} \,\| g\|^2_{{\rm Lip}(f)} ,$$
which is equivalent to (\ref{eqVA}).~$\Box$\bigskip

By Lemma \ref{lemV} and (\ref{eqLIPA}), we can now easily
prove Corollaries \ref{corSD}, \ref{corDC} and \ref{corEA}:
\smallskip

{\bf Proof of Corollary \ref{corSD}.} \ Existence and uniqueness of a stationary distribution $\mu$ for
$(p_t)_{t\ge 0}$ satisfying $\int |y|\, \mu (dy)<\infty $
follows easily as in \cite{KW}, Section~3: By Corollary
\ref{cor16}, the map $\nu\mapsto\nu p_1$ is a contraction
w.r.t.\ the distance $W_f$ (equivalent to $W^1$) on the
complete metric space $\mathcal P^1$ of all probability
measures $\nu$ on $(\rd ,\mathcal B(\rd ))$ satisfying
$\int |y|\, \mu (dy)<\infty $. Hence by the Banach
fixed point theorem, there exists a unique probability
measure $\mu_0$ such that $\mu_0 p_1=\mu_0 $. It is then
elementary to verify that the measure $\mu =\int_0^1
\mu_0p_s\, ds$ satisfies $\mu p_t=\mu $ for any $t\in [0,1]$, and hence for any $t\in [0,\infty )$. Moreover,
by Corollary \ref{cor16},
$$W_f(\mu ,\nu p_t)\ =\ W_f(\mu p_t,\nu p_t)\ \le\
e^{-ct}\, W_f(\mu ,\nu )$$
for any $\nu\in\mathcal P^1$. In particular, as $t\to\infty$, $p_t(x,\cdot )\to\mu$ in $\mathcal P^1$ for
any $x\in\rd$. The variance bound for $\mu$ now follows from the corresponding bound for $p_t(x,\cdot )$ in
Lemma \ref{lemV}.~$\Box$\bigskip

{\bf Proof of Corollary \ref{corDC}.} \ By Lemma \ref{lemV},
\begin{eqnarray*}
\lefteqn{{\rm Cov}\, \left( g(X_t),h(X_{t+s})\right)
\ =\ \mathbb E\left[ g(X_t)\,h(X_{t+s})\right]\, -\,
 E\left[ g(X_t)\right]\,  E\left[ h(X_{t+s})\right] }\\
 &=&\mathbb E\left[ (g\,p_sh)(X_t)\right]\, -\,
 E\left[ g(X_t)\right]\,  E\left[ (p_sh)(X_{t})\right]\
 =\ {\rm Cov}_{p_t(x_0,\cdot )}(g,p_sh)\\ & \le &
 {(1-\exp ({-2ct}))}\, {(2c)^{-1}} \,\| g\|_{{\rm Lip}(f)}\| p_sh\|_{{\rm Lip}(f)}
\end{eqnarray*}
for any $s,t\ge 0$. The assertion now follows by
(\ref{eqLIPA}).~$\Box$\bigskip

{\bf Proof of Corollary \ref{corEA}.} \ The bound for the
bias follows immediately from (\ref{eqLIPA}), since
\begin{eqnarray*}
 \left|\mathbb E\left[\frac 1t \int_0^t g(X_s)\, ds\, -
 \, \int g\, d\mu \right]\right|
&=&\left|\frac 1t \int_0^t \int (p_sg(x_0)-p_sg (y))\,
\mu (dy)\, ds\right|\\
&\le & \frac 1t \int_0^t e^{-cs}\, ds\ \| g\|_{{\rm Lip}(f)}\,  \int d_f(x_0,y)\,
\mu (dy).
\end{eqnarray*}
Moreover, by Corollary \ref{corDC},
\begin{eqnarray*}
 {\rm Var}\, \left(\frac 1t \int_0^t g(X_s)\, ds \right)
&=&{\rm Cov}\, \left(\frac 1t \int_0^t g(X_s)\, ds\, ,
\, \frac 1t \int_0^t g(X_s)\, ds \right)\\
&=&\frac{2}{t^2}\int_0^t\int_s^t{\rm Cov}\,\left(g(X_s),g(X_u)\right)
\, du\, ds\\
&\le & \frac{1}{ct^2}\int_0^t(1-e^{-2cs})\int_s^te^{-c (u-s)}\, du\, ds\; \| g\|^2_{{\rm Lip}(f)}\\
&\le &\frac{1}{c^2t}\, \| g\|^2_{{\rm Lip}(f)}.
\qquad\qquad\qquad\qquad\Box
\end{eqnarray*}

\section{Examples}\label{secexamples}

We now prove the results in Sections
\ref{sub:firstexamples} and \ref{sub:highdim},
including in particular Lemma
\ref{lemP}, Lemma \ref{lemP1} and Theorem \ref{thmLC}.\medskip

{\bf Proof of Lemma \ref{lemP} and Remark \ref{remP}.}
We first prove the lower bounds on the exponential decay rate $c$ in (\ref{eq:10}) stated in (\ref{eq18c}), (\ref{eq18a})
and (\ref{eq18b}). Notice that the constant $c$ defined by
(\ref{eq:10}) increases if $\kappa(r)$ is replaced by a greater function. Indeed, for $r\geq 0$,
\begin{equation}
\label{eq:16}
\Phi(r)\varphi(r)^{-1}\:=\:\int\limits_0^r\varphi(t)\varphi(r)^{-1}\:dt\:=\:\int\limits_0^r\exp\left(\frac 14\int\limits_t^rs\kappa(s)^-\:ds\right)\:dt,
\end{equation}
whence $R_0$, $R_1$ and $c^{-1}=\alpha\int\limits_0^{R_1}\Phi(s)\varphi(s)^{-1}\:ds$ are decreasing functions of $\kappa$.
\medskip\\
{\em Convex Case.} Suppose first that $\kappa(r)\geq 0$ for any $r\geq 0$ and $\kappa(r)\geq K$ for $r\geq R$ with constants $K \in(0,\infty)$ and $R\in[0,\infty)$. Then $R_0=0$, $R_1\leq\max(R,\sqrt{8/K})$, $\varphi\equiv 1$, and hence
\[c\:=\:(\alpha R_1^2/2)^{-1}\:\geq \:\alpha^{-1}\min(R^{-2}/2,K/4).\]
%In particular, if $\kappa(r)\geq K>0$ holds for all $r\geq 0$ then we obtain $c\ge K/(4\alpha)$. This bound for the exponential decay rate is fairly sharp, since in the situation of Example \ref{exlangevin} above with $d=1$ and $U(x)=Kx^2/2$, the exact decay rate is $K/2$. Similarly, if $\kappa(r)=0$ for $r<R$, the order $O(R^{-2})$ for $c$ is optimal as the example of a one-dimensional diffusion with $\sigma=1$ and vanishing drift on $(-R/2,R/2)$ shows.
{\em Locally non-convex case.} Now suppose that $\kappa(r)\geq -L$ for $r\leq R$ and $\kappa(r)\ge K$ for $r>R$ with constants $K,L\in(0,\infty)$ and $R\in[0,\infty]$. Since $\varphi(r)=\varphi(R_0)$ and $\Phi(r)=\Phi(R_0)+(r-R_0)\varphi(R_0)$ for $r\geq R_0$, we have
\begin{eqnarray}
\nonumber
\alpha^{-1}c^{-1} &=&
\int\limits_0^{R_1}\Phi(s)\varphi(s)^{-1}\:ds\\
&=&\int\limits_0^{R_0}\Phi(s)\varphi(s)^{-1}\:ds+(R_1-R_0)\Phi(R_0)\varphi(R_0)^{-1}+(R_1-R_0)^2/2.
\label{eq:p17}\end{eqnarray}
The lower curvature bounds imply the upper bounds
\begin{eqnarray}
R_0&\leq &R,\qquad
R_1-R_0\ \leq\ \min(8/(KR_0),\sqrt{8/K}),\qquad\mbox{and}\label{eq188}\\
\Phi(r)\varphi(r)^{-1}&\leq &\int\limits_0^r\exp(L(r^2-t^2)/8)\:dt\nonumber\\
 &\leq & \min(\sqrt{2\pi/L},r)\exp(Lr^2/8)\qquad\mbox{ for }r\leq R_0.\label{eq189}\end{eqnarray}
Since $\exp x\le 1+(e-1)x$ for $x\in [0,1]$ and
$$\int_0^x\exp (u^2)\, du\ \le\ e+\int_1^x(2-u^{-2})\exp (u^2)\, du\ =\ x^{-1}\exp (x^2)\qquad\mbox{for }x \ge 1,$$
we can conclude that
\begin{eqnarray*}
\int\limits_0^{R_0}\Phi(r)\varphi(r)^{-1}\:dr &\leq &
\int_0^{R_0}r\exp (Lr^2/8)\, dr\ =\ 4L^{-1}
(\exp (LR_0^2/8 )-1)\\ &\le &
(e-1)R_0^2/2\qquad\mbox{ if }LR_0^2/8\leq 1,\qquad\mbox{ \ and}\\
\int\limits_0^{R_0}\Phi(r)\varphi(r)^{-1}\:dr &\leq &
\sqrt{\frac{2\pi}{L} }\int_0^{R_0}\exp (\frac{Lr^2}8)\, dr\ =\
\sqrt{\frac{8\cdot 2\pi}{L^2} }\int_0^{\sqrt{LR_0^2/8}}\exp (u^2)\, du\\
&\le &
8\sqrt{2\pi}L^{-3/2}R_0^{-1}\exp(LR_0^2/8)\qquad\mbox{ \ if }LR_0^2/8\geq1.
\end{eqnarray*}
Combining these estimates, we obtain by (\ref{eq:p17}),
(\ref{eq188}) and (\ref{eq189}),
\begin{eqnarray*}
\alpha^{-1}c^{-1}&\leq &(e-1)R^2/2+e\sqrt{8/K}R+4/K\qquad\mbox{ if }LR_0^2/8\:\leq\:1,\mbox{ \ and }\\
\alpha^{-1}c^{-1}&\leq&8\sqrt{2\pi}R^{-1}L^{-1/2}(L^{-1}+K^{-1})\exp(LR^2/8)+32R^{-2}K^{-2}\mbox{ if }LR_0^2/8\geq 1,
\end{eqnarray*}
where we have used that the function $x\mapsto x^{-1}
\exp (x^2)$ is increasing for $x\ge 1$.~$\Box
$\bigskip

{\bf Proofs for Example \ref{exdoublewell}.}
Consider the one-dimensional Langevin diffusion $(X_t)$
with drift $-\nabla U(x)/2$ and generator
\begin{equation}\label{eqGEN}
\mathcal Lv\ =\ \frac 12(v''-U'v')\ =\ \frac 12\,e^U
\left( e^{-U}v^\prime\right)^\prime .
\end{equation}
The assumption $\liminf_{|x|\to\infty}U^{\prime\prime}(x)>0$ implies that there is a unique strictly positive
bounded eigenfunction $v_1\in C^2(0,\infty )\cap
C([0,\infty ))$ satisfying $v_1(0)=0$, $v_1^\prime (0)=1$
and $\mathcal Lv_1=-\lambda_1v_1$, where
$$\lambda_1\ =\  \lambda_1 (0,\infty )\ =\
\inf_{v\in C_0^\infty(0,\infty )}\frac{\frac 12
\int_0^\infty v'(x)^2\,\exp (-U(x))\, dx}{\int_0^\infty v(x)^2\,\exp (-U(x))\, dx} $$
is the infimum of the spectrum of the self-adjoint
realization of $-\mathcal L$ with Dirichlet boundary
conditions on $(0,\infty )$. Since $\mathcal Lv_1=-
\lambda_1v_1$ and $v_1$ is bounded, the process
$M_t=\exp (\lambda_1t)v_1(X_t)$ is a martingale.
Optional stopping applied to the diffusion with initial
condition $X_0=x_0$ shows that
\begin{eqnarray}
\nonumber
v_1(x_0)& =& \mathbb E_{x_0}\left[ M_0\right]\ =\
\mathbb E_{x_0}\left[ M_{\tau_0\wedge t}\right]\ =\
\mathbb E_{x_0}\left[\exp (\lambda_1t)v_1(X_t);\, {\tau_0> t}\right]\\
\label{eqExbdd}
&\le & \exp (\lambda_1t)\, \mathbb P_{x_0}\left[ {\tau_0> t}\right] \, \sup v_1
\end{eqnarray}
for any $x_0>0$ and $t\ge 0$. Since $v_1(x_0)>0$ and
$\sup v_1 <\infty $, the estimate (\ref{eqExbdd}) implies
the asymptotic lower bound
\begin{equation}
\liminf_{t\to\infty}t^{-1}\log\:P_{x_0}[\tau_0>t]\ \ge\ -\lambda_1(0,\infty)\label{eq:starb} .
\end{equation}
Moreover,
for any fixed $t\le \lambda_1^{-1}/4$,
$$\mathbb P_R \left[ {\tau_0> t}\right]\ \ge\
e^{-1/4}\, v_1(R)/\sup v_1\ \ge \ 3/4 $$
provided $v_1(R)\ge \frac 34e^{1/4}\sup v_1=0.96\ldots
\cdot\sup v_1$. By the eigenfunction equation $e^U(e^{-U}v_1^\prime )'=-\lambda_1v_1$, one verifies that the
latter condition is satisfied whenever $U$ is growing fast
enough on $[R,\infty )$.
%XXXCheck
\medskip

For bounding $\lambda_1(0,\infty )$ from above let
$$v(x)\ =\ \min (\sqrt Lx,1)\ =\ \left\{
\begin{array}{ll}
\sqrt L x &\mbox{ if }x\le 1/\sqrt L,\\
1 &\mbox{ if }x\ge 1/\sqrt L.
\end{array}\right.$$
By the assumptions on $U$, the function $v$ is contained
in the weighted Sobolev space $H_0^{1,2}((0,\infty ),
e^{-U}\, dx)$ (closure of $C_0^\infty (0,\infty )$ w.r.t.\
the norm $\| w\|^2=\int_0^\infty (w^2+(w')^2)\, e^{-U}\, dx$). Therefore, if $LR^2/4\ge 1$ then
(\ref{eq:star}) holds, since
%
%
%
%
%
%
%
%Consider the Langevin diffusion case in $\mathbb R^1$ a symmetric potential $U(x)$ satisfying $U(x)=-Lx^2/2$ for $x\in[-R/2,R/2]$. If $\|\:\cdot\:\|$ is the Euclidean norm then by (\ref{eq:14}), $\kappa(r)\leq -L$ for $r\in(0,R]$. Moreover, for an arbitrary concave strictly increasing function $f\in C(\mathbb R^+)$, the $L^1$ Wasserstein distance $W_f(\delta_{-R/2}p_t,\delta_{R/2}p_t)$ between the distributions at time $t$ of the diffusion starting at $-R/2$ and $R/2$ is at least of order $P_{R/2}[T_0>t]$, where $T_0$ denotes the first hitting time of $0$ for the process starting at $R/2$. By a well-known large deviation principle [\ref{???}],
%\[\lim\limits_{t\to\infty}t^{-1}\log\:P_{R/2}[T_0>t]\:=\:-\lambda_1(0,\infty),\]
%where $\lambda_1(0,\infty)$ denotes the lowest Dirichlet eigenvalue of the generator on $(0,\infty)$. Inserting the function $g(x)=\min(\sqrt Lx,1)$ in the variational characterization of the Dirichlet eigenvalue shows that
\begin{eqnarray*}
\lambda_1 &\leq &\frac{\frac 12\int v'(x)^2\exp(-U(x))\:dx}{\int v(x)^2\exp(-U(x))\:dx}\ \le\
\frac{\int_0^{1/\sqrt L}L\exp(Lx^2/2)\,dx}{\int_0^{R/2}v(x)^2\exp(Lx^2/2)\,dx}\\
&= &\frac L2\,\frac{\int_0^{1}exp(y^2/2)\,dy}{\int_0^{\sqrt{LR^2/4}}\min (y,1)^2\exp(y^2/2)\,dy}
\  \le\
\frac{3Le^{1/2}}{2}\, \sqrt{\frac{LR^2}{4}}\, \exp\left(\frac{LR^2}{8} \right) .
\end{eqnarray*}
Here we have used that by assumption, $U(x)\ge -Lx^2/2$
for any $x\in\mathbb R$ with equality for $|x|<R/2$, and for $x\ge 1$,
$$\int_0^x\min (y,1)^2e^{y^2/2}\, dy\ =\ \int_0^1\ldots
+\int_1^x\ldots\ \ge\ \frac 13 +\frac 1xe^{x^2/2}-1
\ \ge\ \frac{1}{3x}e^{x^2/2}$$
as $(x^{-1}e^{x^2/2})'=(1-x^{-2})
e^{x^2/2}\le e^{x^2/2}$.$~\Box
$\bigskip

{\bf Proof of Lemma \ref{lemP1}.} Since $b=b_0+\gamma$,
we have
$$(x-y)\cdot G(b(x)-b(y))\ =\
(x-y)\cdot G(b_0(x)-b_0(y))+(x-y)\cdot G(\gamma (x)-\gamma (y))$$
for any $x,y\in\rd $. Therefore, by (\ref{eqPC}) and by
definition of $\kappa$ and $\kappa_0$,
\begin{eqnarray}
\label{eqPE} \kappa (r)^-&\le & \kappa_0(r)^-
\qquad\qquad\mbox{ for any }r\le R,\mbox{ and}\\
\label{eqPD} \kappa (r)^-&\le & \kappa_0(r)^-
+4r^{-1}\sup\|\gamma\|\qquad\mbox{for any }r\in (0,\infty ).
\end{eqnarray}
In particular, if $\gamma $ is bounded then $\kappa $
satisfies the conditions in (\ref{eqassk}). Since the
constant $R_1(b)$ defined w.r.t.\ $b$ is smaller than
the corresponding constant $R_1$ defined w.r.t.\ $b_0$,
we obtain
\begin{eqnarray*}
\frac 1c &\le & \int_0^{R_1}\int_0^s\exp\left(
\frac 14\int_t^su\kappa (u)^-\, du\right)\, dt\, ds\\
&\le & \int_0^{R_1}\int_0^s\exp\left(
\frac 14\int_t^su\kappa_0 (u)^-\, du\right)\,
\exp\left( R\sup\|\gamma\|\right) \, dt\, ds\\
&\le & \frac{1}{c_0}\cdot \exp\left( R\sup\|\gamma\|\right) ,
\end{eqnarray*}
i.e., (\ref{eqPF}) holds.\medskip

Similarly, if $\gamma $ satisfies the one-sided Lipschitz
condition (\ref{eqPG}) then
\begin{equation}
\label{eqPI}
\kappa (r)^-\ \le \ \kappa_0(r)^-
+2L\qquad\mbox{for any }r\in (0,\infty ).
\end{equation}
Hence again the conditions in (\ref{eqassk}) are satisfied,
and we obtain
$$\frac 1c\ \le \ \frac{1}{c_0}\cdot \exp\left( \frac L2\int_0^Rr\, dr\right) $$
similarly as above, i.e., (\ref{eqPH}) holds.~$\Box
$\medskip

{\bf Proof of Theorem \ref{thmLC}.}
Fix $R>0$ and probability measures $\mu ,\nu$ on $\rd$.
By definition of $f_R$,
$$f_R''(r)\ \le\ \frac 14 r\kappa (r)f_R'(r)-f_R(r)
\left/ \int_0^R\frac{\Phi (s)}{\varphi (s)}\, ds\right.$$
for any $r<R$. Therefore, similarly to the proof of
Theorem \ref{thmmain1}, Equation (\ref{eq:13}) shows
that the process $e^{c_Rt}f_R(r_t)$ is a local supermartingale for $t<\hat\tau_R$ where
$$\hat\tau_R\ =\ \inf\{ t\ge 0: r_t> R\} .$$
Here $r_t=\| X_t-Y_t\| $ again denotes the distance process
for a reflection coupling $(X_t,Y_t)$ of two solutions of
(\ref{eq:1}) with initial distribution given by a coupling
$\eta$ of $\mu $ and $\nu$. By optional stopping and Fatou's
lemma, we thus obtain
$$
\mathbb{E}[f_R(r_t);\, \hat\tau_R >t]\ \le\
\mathbb{E}[f_R(r_{t\wedge\hat\tau_R})]\ \le\
\exp ({-c_Rt})\,\mathbb{E}[f_R(r_0)]
$$
for any $t\ge 0$, and hence
\begin{eqnarray*}
\mathbb{E}[f_R(r_t)]&\le & \exp ({-c_Rt})\mathbb{E}[f_R(r_0)]
\,+\,\mathbb P [ \hat\tau_R \le t]\\ & \le &
e^{-c_Rt}\int f_R(\| x-y\|\, \eta (dx\,dy)\, +\,
\mathbb P_\mu [\tau_{R/2}\le t]\, +\,
\mathbb P_\nu [\tau_{R/2}\le t] .
\end{eqnarray*}
The assertion now follows as in the proof of Corollary
\ref{cor16} by minimizing over all couplings $\eta$
of $\mu$ and $\nu$.~$\Box$\medskip

\section{Couplings on product spaces}\label{secprod}

Let $d=\sum_{i=1}^nd_i$ with $n,d_1,\ldots ,d_n\in\mathbb{N}$. We now consider ``componentwise''
couplings for diffusion processes $X_t=(X_t^1,\ldots ,
X_t^n)$ and $Y_t=(Y_t^1,\ldots ,
Y_t^n)$ on $\rd $ satisfying the s.d.e.
\begin{equation}
\label{eqH}
dX_t^i\ =\ b^i(X_t)\, dt\, +\, dB_t^i,\qquad i=1,\ldots ,n,
\end{equation}
with initial conditions $X_0\sim\mu $ and $Y_0\sim\nu $.
Here $B^i$, $i=1,\ldots ,n$, are independent Brownian
motions on $\mathbb R^{d_i}$, and $b^i:\mathbb R^{d_i}\to
\mathbb R^{d_i}$ are locally Lipschitz continuous functions
such that the unique strong solution of (\ref{eqH}) is
non-explosive for any given initial condition.\smallskip

Let $\delta >0$. Suppose that $\lambda^i,\pi^i:\rd
\to [0,1]$, $i=1,\ldots ,n$, are Lipschitz continuous
functions such that
\begin{eqnarray}
\label{eqI}
\lambda^i(z)^2+\pi^i(z)^2 &=& 1\qquad\mbox{for any }z\in\rd
,\qquad\mbox{and}\\
\label{eqJ}
\lambda^i(z) &=& 0\qquad\mbox{if }|z^i|\le\delta /2,
\end{eqnarray}
and let $B^i$ and $\tilde B^i$, $1\le i\le n$, be independent Brownian motions on $\mathbb R^{d_i}$. Then
a coupling of two solutions of (\ref{eqH}) with initial
distributions $\mu$ and $\nu$ is given by a strong solution
of the system
\begin{eqnarray}
\label{eqK}
dX_t^i & =& b^i(X_t)\, dt\, +\,\lambda^i(Z_t)\,  dB_t^i\,
+\,\pi^i(Z_t)\, d\tilde B_t^i ,\\
\nonumber
dY_t^i & =& b^i(Y_t)\, dt\, +\,\lambda^i(Z_t)\,(I-2e_t^i
e_t^{i,T})\,   dB_t^i\,
+\,\pi^i(Z_t)\, d\tilde B_t^i ,
\end{eqnarray}
$1\le i\le n$, with initial distribution $(X_0,Y_0)\sim
\eta$ where $\eta$ is a coupling of $\mu$ and $\nu$. Here
we use the notation $$Z_t\ =\ X_t-Y_t,$$
and $e_t^i$ is a measurable process taking values in the
unit sphere in $\mathbb{R}^{d_i}$ such that
$$e_t^i\ =\
\left\{\begin{array}{ll}Z_t^i/|Z_t^i| &\ \mbox{ if }
Z_t^i\neq 0,\\
u^i&\ \mbox{ if }
Z_t^i= 0,
\end{array}\right.$$
where $u^i$ is an arbitrary fixed unit vector in $\mathbb{R}^{d_i}$. Notice that by (\ref{eqJ}), the choice
of $u^i$ is not relevant for (\ref{eqK}), which is a standard It\^o s.d.e.\ in $\mathbb{R}^{2d}$ with locally
Lipschitz continuous coefficients. To see that (\ref{eqK})
defines a coupling, we observe that $(X_t)$ and $(Y_t)$
satisfy (\ref{eqH}) w.r.t.\ the processes
$\hat B_t=(\hat B_t^1,\ldots , \hat B_t^n)$ and
$\check B_t=(\check B_t^1,\ldots , \check B_t^n)$
defined by
\begin{eqnarray*}
\hat B_t^i &=&  \int_0^t\lambda^i(Z_s)\, dB_s^i\,+\,
 \int_0^t\pi^i(Z_s)\, d\tilde B_s^i,\\
\check B_t^i &=&  \int_0^t\lambda^i(Z_s)\, (I-2e_s^i
e_s^{i,T})\, dB_s^i\,+\,
 \int_0^t\pi^i(Z_s)\, d\tilde B_s^i.
\end{eqnarray*}
By L\'evy's characterization and by (\ref{eqI}), both
$\hat B$ and $\check B$ are indeed Brownian motions in $\rd$, cp.\ the corresponding argument for reflection coupling.

\begin{remark}
(1) By Condition (\ref{eqJ}) and non-explosiveness
of (\ref{eqH}), the coupling process $(X_t,Y_t )$ is
defined for any $t\ge 0$.\smallskip\\
(2) By choosing $\lambda^i\equiv 0$ and $\pi^i\equiv 1$
we recover the {\em synchronuous coupling}, i.e., the same noise
is applied to both processes $X$ and $Y$.\smallskip\\
(3) A {\em componentwise reflection coupling} would be
informally given by choosing $\lambda^i(z)=1$ if $z^i\neq 0$ and $\lambda^i(z)=0$ if $z^i=0$. As this function is not continuous and $e^i(z)=z^i/|z^i|$ also has a discontinuity at zero, it is not obvious how to make sense of this
coupling rigorously. Instead, we will use below an
approximate componentwise reflection coupling where
$\lambda^i(z)=1$ if $|z^i|\ge\delta $ and
$\lambda^i(z)=0$ if $|z^i|\le\delta /2 $
for a small positive constant $\delta $.
\end{remark}

By subtracting the equations for $X$ and $Y$ in (\ref{eqK}), we see that the difference process
$Z=X-Y$ satisfies the s.d.e.
\begin{equation}
\label{eqM}
dZ_t^i \ =\ (b^i(X_t)-b^i(Y_t))\, dt\, +\,2\lambda^i(Z_t)\,  e_t^i\, dW_t^i\, ,
\end{equation}
$i=1,\ldots ,n$, where the processes
$$W_t^i\ =\ \int_0^te_t^{i}\cdot dB_t^i,\qquad 1\le i\le n,$$
are independent {\em one-dimensional} Brownian motions.
\smallskip

Let $r_t^i=|X_t^i-Y_t^i|$ denote the Euclidean norm of
$Z_t^i$. The next lemma is crucial for quantifying
contraction properties of the coupling given by (\ref{eqK}):

\begin{lemma}\label{lem3A}
Suppose that $f:[0,\infty )\to [0,\infty )$ is a strictly
increasing concave function in $C^1([0,\infty ))$ such
that $f'$ is absolutely continuous on $(0,\infty )$.
Then for any $i=1,\ldots ,n$, the process $f(r_t^i)$
satisfies the It\^o equation
\begin{eqnarray}
\nonumber f(r_t^i) &=& f(r_0^i)\, +\, 2\int_0^t\lambda^i
(X_s-Y_s)\, f'(r_s^i)\, dW_s^i\\
&&+\int_0^t\left\{e_s^i\cdot (b^i(X_s)-b^i(Y_s))\, f'(r_s^i)\, +\,2\lambda^i(X_s-Y_s)^2\, f''(r_s^i)\right\}\, ds.\label{eqNA}
\end{eqnarray}
\end{lemma}

\begin{remark}
The lemma shows in particular that the process $r_t^i$
satisfies
\begin{equation}
\label{eqN}
dr_t^i\ =\ e_t^i\cdot (b^i(X_t)-b^i(Y_t))\, dt\, +\,2\lambda^i(X_t-Y_t)\, dW_t^i .
\end{equation}
Notice that in this equation, the drift term does not
depend on the choice of $\lambda $.
\end{remark}

{\bf Proof of Lemma \ref{lem3A}.} \ Recall that $e_t^i=
Z_t^i/|Z_t^i|$ if $r_t^i=|Z_t^i|\neq 0$. Since the
function $y\mapsto y/|y|$ is smooth on $\mathbb{R}^{d_i}
\setminus\{ 0\} $ and $x\mapsto\sqrt x$ is smooth on
$(0,\infty )$, we can apply It\^o's formula and (\ref{eqM})
to show that the It\^o equations
\begin{eqnarray}
\nonumber d|Z^i|^2&  =& 2Z^i\cdot (b^i(X)-b^i(Y))\, dt\, +\, 4\, \lambda^i(Z)^2\, dt\, +\,4\lambda^i(Z)\,|Z^i|  \, dW^i\, ,\\
\nonumber
dr^i &=& \frac{1}{2r^i}\, d|Z^i|^2\, - \,\frac{1}{8(r^i)^3}
\, d[|Z^i|^2]\\
\label{eqO}
&=& e^i\cdot (b^i(X)-b^i(Y))\, dt\,  +\,2\lambda^i(X-Y)\, dW^i
\end{eqnarray}
hold almost surely on any stochastic interval $[\tau_1,\tau_2]$ such that $Z_t^i\neq 0$ a.s.\ for $\tau_1\le t
\le\tau_2$.\smallskip

On the other hand, suppose that $|Z^i| <\delta /2$ a.s.\ on
a stochastic interval $[\tau_3,\tau_4 ]$. Then on
$[\tau_3,\tau_4 ]$, $\lambdaî (Z)\equiv 0$ by (\ref{eqJ}),
and hence $Z^i$ is almost surely absolutely continuous
with
$$dZ^i/dt\ =\ b^i(X)-b^i(Y)\qquad\mbox{a.e.\ on}\ [\tau_3,\tau_4 ].$$
This implies that $r^i=|Z^i|$ is almost surely absolutely continuous on $[\tau_3,\tau_4 ]$ as well with
\begin{equation}
\label{eqP}
dr^i/dt\ =\ e^i\cdot (b^i(X)-b^i(Y))\qquad\mbox{a.e.\ on}
\  [\tau_3,\tau_4 ],
\end{equation}
which is equivalent to (\ref{eqN}) on $[\tau_3,\tau_4 ]$.
Note that the value of $e^i$ for $Z^i=0$ is not relevant here, since $Z^i$ can only stay at $0$ for a positive
amount of time if $b^i(X)-b^i(Y)$ vanishes during that
time interval.\smallskip

Since $\mathbb{R}_+$ is the union of countably many stochastic intervals of the first and second type
considered above, the It\^o equation (\ref{eqN}) holds
almost surely on $\mathbb{R}_+$. The assertion
(\ref{eqNA}) now follows from (\ref{eqN}) by another
application of It\^o's formula. Here it is enough to assume
that $f$ is $C^1$ on $[0,\infty )$ and $f'$ is absolutely
continuous on $(0,\infty )$ because $\lambda^i(X_s-Y_s)$
vanishes for $r_s^i<\delta /2$.~$\Box$\medskip

We now fix weights $w_1,\ldots w_n\in [0,\infty )$ and
strictly increasing concave functions $f_1,\ldots ,f_n
\in C^1 ([0,\infty ))\cap C^2((0,\infty ))$ such that
$f_i(0)=0$ for any $i$. Consider
\begin{equation}
\label{eqQ}
\rho_t\ =\ \sum_{i=1}^nf_i(r_t^i)\, w_i\ =\ d_{f,w}(X_t,Y_t)
\end{equation}
where $d_{f,w}$ is defined by (\ref{eq14delta}). By
Lemma \ref{lem3A},
\begin{eqnarray}
\nonumber d\rho_t &=& \sum_{i=1}^n \left( e_t^i\cdot (b^i(X_t)-b^i(Y_t))\, f_i^\prime (r_t^i)\,+\,
 2\lambda^i(X_t-Y_t)^2\, f_i^{\prime\prime}(r_t^i)\right)
 \, w_i\, dt\\
&&+ 2\,\sum_{i=1}^n\lambda^i
(X_t-Y_t)\, f_i^\prime (r_t^i)\, dW_t^i .
\label{eqR}
\end{eqnarray}
Notice that the last term on the right hand side is a
martingale since $\lambda^i$ and $f_i^\prime $ are bounded.
This enables us to control the expectation $\mathbb E[\rho_t]$ if we can bound the drift
in (\ref{eqR}) by $m-c\rho_t$ for constants
$m,c\in (0,\infty )$:

\begin{lemma}\label{lem3C}
Let $m,c\in (0,\infty )$ and suppose that
\begin{equation}
\label{eqQB}
\sum_{i=1}^n \left( cf_i(r^i)\, +\, (x^i-y^i)\cdot (b^i(x)-b^i(y))\, \frac{f_i^\prime (r^i)}{r^i} \,+\,
 2\lambda^i(x-y)^2\, f_i^{\prime\prime}(r^i)\right)
 \, w_i\ \le\ m
\end{equation}
holds for any $x,y\in\rd $ with $r^i:=|x^i-y^i|>0\ \forall\,
i\in\{ 1,\ldots n\} $. Then
\begin{equation}
\label{eqS}
\mathbb{E}[\rho_t]\ \le\ e^{-ct}\,\mathbb{E}[\rho_0]
\, +\, m\, (1-e^{-ct})/c\qquad\mbox{for any }t\ge 0.
\end{equation}
\end{lemma}

{\bf Proof.} \ We first note that by continuity of $b^i$
and $f_i^\prime$, (\ref{eqQB}) implies that
\begin{equation}
\label{eqQA}
\sum_{i=1}^n \left( cf_i(r^i)\, +\, e^i\cdot (b^i(x)-b^i(y))\, {f_i^\prime (r^i)} \,+\,
 2\lambda^i(x-y)^2\, f_i^{\prime\prime}(r^i)\right)
 \, w_i\ \le\ m
\end{equation}
holds for any $x,y\in\rd $ (even if $x^i-y^i=0$) provided
$e^i=(x^i-y^i)/r^i$ if $r^i>0$ and $e^i$ is an arbitrary
unit vector if $r^i=0$. Indeed, we obtain (\ref{eqQA})
by applying (\ref{eqQB}) with $x^i$ replaced by $x^i+he^i$
whenever $x^i-y^i=0$ and taking the limit as $h\downarrow 0$. In particular, by (\ref{eqQA}), the drift term $\beta_t$ in (\ref{eqR})
is bounded from above by
$$\beta_t\ \le\ m-\sum_{i=1}^ncf_i(r_t^i)w_i\ =\ m-c\rho_t .$$
Therefore by (\ref{eqR}) and by the It\^o product rule,
$$d(e^{ct}\rho )\ =\ e^{ct}\, d\rho\, +\,
ce^{ct}\rho\, dt\ \le\ e^{ct}m\, dt\, +\, dM$$
where $M$ is a martingale, and thus
$$\mathbb{E}[e^{ct}\rho_t]\ \le\ \mathbb{E}[\rho_0]
\, +\, m\, \int_0^te^{cs}\, ds\qquad\mbox{for any }t\ge 0.
$$
\begin{flushright}
$\Box$
\end{flushright}
\medskip

Since $f_i^{\prime\prime}\le 0$, the process $\rho_t$ is
decreasing more rapidly (or growing more slowly) if $\lambda^i$ takes larger values. In particular, the decay
properties of $\rho_t$ would be optimized when $\lambda^i(z)=1 $ for any $z$ with $z^i\neq 0$. This
optimal choice of $\lambda^1,\ldots ,\lambda^n$ would
correspond to a componentwise reflection coupling, but
it violates Condition (\ref{eqJ}). It is perhaps possible
to construct a corresponding coupling process by an
approximation argument. For our purpose of bounding the
Kantorovich distance $W_{f,w}( \mu p_t,\nu p_t)$ this
is not necessary. Indeed, it will be sufficient to
consider approximate componentwise reflection couplings
where (\ref{eqI}) and (\ref{eqJ}) are satisfied and
$\lambda^i(z)=1$ whenever $|z^i|>\delta$. The limit
$\delta\downarrow 0$ will then be considered for the
resulting estimates of the Kantorovich distance but not
for the coupling processes.

\section{Application to interacting diffusions}\label{secinteract}

We will now apply the couplings introduced in Section
\ref{secprod} to prove the contraction properties for
systems of interacting diffusions stated in Theorem
\ref{thmA} and Corollary \ref{corB}.
We consider the setup described in Section \ref{subsectioncrc}, i.e.,
\begin{equation}
\label{eqSA}b^i(x)\ =\ b_0^i(x^i)\, +\,\gamma^i(x)
\qquad\mbox{for } i=1,\ldots ,n
\end{equation}
with $b_0^i:\mathbb{R}^{d_i}\to\mathbb{R}^{d_i}$ locally
Lipschitz such that $\kappa_i$ defined by (\ref{eqKI}) is
continuous on $(0,\infty )$ with
\begin{equation}\label{eqAK}
\liminf_{r\to\infty}\kappa_i(r)>0\ \mbox{ and }\
\lim_{r\to 0}r\kappa_i(r)=0\ \mbox{ for any }1\le i\le n.
\end{equation}
The functions $f_i$ are defined via $\kappa_i$, and $c_i$
is the corresponding contraction rate given by
(\ref{eq14dot}).\medskip

{\bf Proof of Theorem \ref{thmA}.} \ We fix $\delta >0$ and
Lipschitz continuous functions $\lambda^i,\mu^i:\rd\to
[0,1]$, $1\le i\le n$, such that (\ref{eqI}) and (\ref{eqJ}) hold and $\lambda^i(z)=1$ if $|z^i|\ge\delta $.
Let $(X_t,Y_t)$ denote a corresponding approximate
componentwise reflection coupling of two solutions of
(\ref{eqcircstar}) given by (\ref{eqK}), and let
$\rho_t=d_{f,w}(X_t,Y_t)$. We will apply Lemma \ref{lem3C}
which requires bounding the right hand side in (\ref{eqQB}). For this purpose recall that $f_i$ and $c_i$
have been chosen in such a way that
$$2f_i^{\prime\prime}(r)-\frac 12 r\kappa_i(r)f_i^\prime
(r)\ \le\ -c_i\, f_i(r)\qquad\forall\ r>0,$$
cf.\ (\ref{eq:14}) and (\ref{eq:15}). Therefore,
by (\ref{eqSA}) and by definition of $\kappa_i$,
\begin{eqnarray}
\nonumber
\lefteqn{(x^i-y^i)\cdot (b^i(x)-b^i(y))\, f_i^\prime (r^i)/r^i\, +\, 2\lambda^i(x-y)^2\, f_i^{\prime\prime }(r^i)}\\
\nonumber &\le & -\frac 12r^i\kappa_i(r^i)f_i^\prime (r^i)
\, +\, |\gamma^i(x)-\gamma^i(y)|f_i^\prime (r^i)\, +\,
2\lambda^i(x-y)^2\, f_i^{\prime\prime }(r^i)\\
\label{eqT} &\le & -\lambda^i(x-y)^2 c_i f_i(r^i)
+ |\gamma^i(x)-\gamma^i(y)| -
\frac 12(1-\lambda^i(x-y)^2)\,r^i\kappa_i(r^i) f_i^{\prime }(r^i)\\
\nonumber &\le & - c_if_i(r^i)
\, +\, |\gamma^i(x)-\gamma^i(y)|\, +\,
c_i\delta\, +\, \frac 12\sup_{r<\delta }\left(r\kappa_i(r)^-\right)
\end{eqnarray}
for any $x,y\in\rd$ with $r^i=|x^i-y^i|>0$. Here we have
used that $0\le f_i^\prime \le 1 $, and that $\lambda^i
(x-y)\neq 1 $ only if $r^i<\delta $. In this case,
$f_i(r^i)\le r^i\le\delta $. By (\ref{eqT}) and by the
assumption (\ref{eq14dotdot}) on $\gamma^i$, we obtain
\begin{eqnarray*}
\lefteqn{\sum_{i=1}^n\left( (x^i-y^i)\cdot (b^i(x)-b^i(y))\, f_i^\prime (r^i)/r^i\, +\, 2\lambda^i(x-y)^2\, f_i^{\prime\prime }(r^i)\right)\, w_i}\\
&\le & m(\delta )\, +\, \sum_{i=1}^n(-c_i+\varepsilon_i)
f_i(r^i)w_i\ \le\ m(\delta )\, -\, c\sum_{i=1}^nf_i(r^i)w_i
\end{eqnarray*}
for $x,y$ as above, where
$$m( \delta )\ =\ \sum_{i=1}^n (c_i\delta +\frac 12\sup_{r<\delta }(r\kappa_i(r)^-)$$
is a finite constant by (\ref{eqAK}), and $c=\min_{i=1,\ldots n}(c_i-\varepsilon_i)$. Hence (\ref{eqQB}) is satisfied with $c$ and $m(\delta )$ and, therefore,
\begin{equation}
\label{eqUA}
\mathbb E[\rho_t]\ \le\ e^{-ct}\, \mathbb E[\rho_0]\,
+\, m(\delta )\, (1-e^{-ct})/c .
\end{equation}
By choosing the coupling process $(X_t,Y_t)$ with initial
distribution given by a coupling $\eta$ of probability
measures $\mu$ and $\nu$ on $\rd$, we conclude that
\begin{eqnarray}
\nonumber
W_{f,w}(\mu p_t,\nu p_t) &\le &\mathbb E\left[d_{f,w}(X_t,Y_t)\right]\ =\ \mathbb E[\rho_t]\\
\label{eqU} &\le & e^{-ct}\int d_{f,w}(x,y)\,\eta (dx~dy)\,
+\, m(\delta )\, (1-e^{-ct})/c
\end{eqnarray}
for any $t\ge 0$. Moreover, by (\ref{eqAK2}),
$m(\delta )\to 0$ as $\delta\downarrow 0$. Hence the assertion (\ref{eqAX}) follows from (\ref{eqU}) by taking
the limit as $\delta\downarrow 0$ and minimizing over
all couplings $\eta$ of $\mu$ and $\nu$. Finally,
(\ref{eqAY}) follows from (\ref{eqAX}) since
$\varphi (R^i_0)r/2\le f_i(r)\le r$ implies
$$A^{-1}\, d_{\ell^1}(x,y)\ \le\ d_{f,w}(x,y)\ =\
\sum f_i(|x^i-y^i|)\, w_i\ \le\ d_{\ell^1}(x,y).\
\ \Box $$

\

{\bf Proof of Corollary \ref{corB}.} \ The $\ell^1 $-Lipschitz condition (\ref{eqA}) for $\gamma$ implies
that (\ref{eq14dotdot}) holds with $w_i=1$ for any $i$,
and
$$\lambda\varepsilon_i^{-1}\ =\ \inf_{r>0}f_i(r)\ =\
f_i^\prime (R_1^i)\ =\ \varphi_i (R_0^i)/2,$$
i.e., $\varepsilon_i=2\lambda /\varphi_i(R_0^i)$.
The assertion now follows from Theorem \ref{thmA}.
\begin{flushright}
$\Box$
\end{flushright}
\medskip

\begin{acknowledgements}
I would like to thank the referees for helpful comments.
Financial support from the German Science Foundation through the
{\em Hausdorff Center for Mathematics} is gratefully acknowledged.
\end{acknowledgements}

\bibliographystyle{spmpsci}

\bibliography{bibcoupling}

%\begin{thebibliography}{00}
%\bibitem{CL}
%{M.F.~Chen, S.F.~Li}, {Coupling methods for multidimensional diffusion processes},
%   {Ann. Probab.}~{17}
%      {(1989)} {151--177}.
%\bibitem{CW}
%    {M.F.~Chen, F.-Y.~Wang},
%      {Estimation of spectral gap for elliptic operators},
%   {Trans. Amer. Math. Soc.}~{349}
%      {(1997)}
%      {1239--1267}.
%\bibitem{F}
%  {M.~Freidlin},
%     {Functional integration and partial differential equations},
% {Princeton University Press},
%   {Princeton}
%     {1985}.
%\bibitem{HM}
%  {M.~Hairer, J.C.~Mattingly},
%     {Spectral gaps in {W}asserstein distances and the 2{D}
%              stochastic {N}avier-{S}tokes equations},
%   {Ann. Probab.}~{36} {(2008)}
%     {2050--2091}.
%\bibitem{LR}
% {T.~Lindvall, L.C.G.~Rogers}, {Coupling of multidimensional diffusions by reflection},
%   {Ann. Probab.}~{14} {(1986)} {860--872}.
%\bibitem{W}
%  {F.-Y.~Wang},
%     {Functional inequalities, {M}arkov processes and spectral theory},
% {Science Press},
%   {Beijing}
%     {2004}.
%\end{thebibliography}
\end{document}